\documentclass[11pt]{article}
\usepackage{amsfonts,mathrsfs,amssymb,amsthm,mathptm}
\usepackage{amsmath,amscd}
\usepackage{mathptm,pslatex}
\usepackage{color}
\usepackage[dvips]{graphicx}
\usepackage[all]{xy}
\usepackage{graphicx}
\usepackage{fancyhdr}

\begin{document}
\newtheorem{Def}{Definition}[section]
\newtheorem{Bsp}[Def]{Example}
\newtheorem{Prop}[Def]{Proposition}
\newtheorem{Theo}[Def]{Theorem}
\newtheorem{Lem}[Def]{Lemma}
\newtheorem{Koro}[Def]{Corollary}
\theoremstyle{definition}
\newtheorem{Rem}[Def]{Remark}

\newcommand{\add}{{\rm add}}
\newcommand{\con}{{\rm con}}
\newcommand{\gd}{{\rm gl.dim}}
\newcommand{\sd}{{\rm st.dim}}
\newcommand{\sr}{{\rm sr}}
\newcommand{\dm}{{\rm dom.dim}}
\newcommand{\cdm}{{\rm codomdim}}
\newcommand{\tdim}{{\rm dim}}
\newcommand{\E}{{\rm E}}
\newcommand{\Mor}{{\rm Morph}}
\newcommand{\End}{{\rm End}}
\newcommand{\ind}{{\rm ind}}
\newcommand{\rsd}{{\rm res.dim}}
\newcommand{\rd} {{\rm rd}}
\newcommand{\ol}{\overline}
\newcommand{\overpr}{$\hfill\square$}
\newcommand{\rad}{{\rm rad}}
\newcommand{\soc}{{\rm soc}}
\renewcommand{\top}{{\rm top}}
\newcommand{\pd}{{\rm pdim}}
\newcommand{\id}{{\rm idim}}
\newcommand{\fld}{{\rm fdim}}
\newcommand{\Fac}{{\rm Fac}}
\newcommand{\Gen}{{\rm Gen}}
\newcommand{\fd} {{\rm fin.dim}}
\newcommand{\Fd} {{\rm Fin.dim}}
\newcommand{\Pf}[1]{{\mathscr P}^{<\infty}(#1)}
\newcommand{\DTr}{{\rm DTr}}
\newcommand{\cpx}[1]{#1^{\bullet}}
\newcommand{\D}[1]{{\mathscr D}(#1)}
\newcommand{\Dz}[1]{{\mathscr D}^+(#1)}
\newcommand{\Df}[1]{{\mathscr D}^-(#1)}
\newcommand{\Db}[1]{{\mathscr D}^b(#1)}
\newcommand{\C}[1]{{\mathscr C}(#1)}
\newcommand{\Cz}[1]{{\mathscr C}^+(#1)}
\newcommand{\Cf}[1]{{\mathscr C}^-(#1)}
\newcommand{\Cb}[1]{{\mathscr C}^b(#1)}
\newcommand{\Dc}[1]{{\mathscr D}^c(#1)}
\newcommand{\K}[1]{{\mathscr K}(#1)}
\newcommand{\Kz}[1]{{\mathscr K}^+(#1)}
\newcommand{\Kf}[1]{{\mathscr  K}^-(#1)}
\newcommand{\Kb}[1]{{\mathscr K}^b(#1)}
\newcommand{\DF}[1]{{\mathscr D}_F(#1)}

\newcommand{\Kac}[1]{{\mathscr K}_{\rm ac}(#1)}
\newcommand{\Keac}[1]{{\mathscr K}_{\mbox{\rm e-ac}}(#1)}

\newcommand{\modcat}{\ensuremath{\mbox{{\rm -mod}}}}
\newcommand{\Modcat}{\ensuremath{\mbox{{\rm -Mod}}}}
\newcommand{\Spec}{{\rm Spec}}

\newcommand{\stmc}[1]{#1\mbox{{\rm -{\underline{mod}}}}}
\newcommand{\Stmc}[1]{#1\mbox{{\rm -{\underline{Mod}}}}}
\newcommand{\prj}[1]{#1\mbox{{\rm -proj}}}
\newcommand{\inj}[1]{#1\mbox{{\rm -inj}}}
\newcommand{\Prj}[1]{#1\mbox{{\rm -Proj}}}
\newcommand{\Inj}[1]{#1\mbox{{\rm -Inj}}}
\newcommand{\PI}[1]{#1\mbox{{\rm -Prinj}}}
\newcommand{\GP}[1]{#1\mbox{{\rm -GProj}}}
\newcommand{\GI}[1]{#1\mbox{{\rm -GInj}}}
\newcommand{\gp}[1]{#1\mbox{{\rm -Gproj}}}
\newcommand{\gi}[1]{#1\mbox{{\rm -Ginj}}}

\newcommand{\opp}{^{\rm op}}
\newcommand{\otimesL}{\otimes^{\rm\mathbb L}}
\newcommand{\rHom}{{\rm\mathbb R}{\rm Hom}\,}
\newcommand{\pdim}{\pd}
\newcommand{\Hom}{{\rm Hom}}
\newcommand{\Coker}{{\rm Coker}}
\newcommand{ \Ker  }{{\rm Ker}}
\newcommand{ \Cone }{{\rm Con}}
\newcommand{ \Img  }{{\rm Im}}
\newcommand{\Ext}{{\rm Ext}}
\newcommand{\StHom}{{\rm \underline{Hom}}}
\newcommand{\StEnd}{{\rm \underline{End}}}

\newcommand{\KK}{I\!\!K}

\newcommand{\gm}{{\rm _{\Gamma_M}}}
\newcommand{\gmr}{{\rm _{\Gamma_M^R}}}

\def\vez{\varepsilon}\def\bz{\bigoplus}  \def\sz {\oplus}
\def\epa{\xrightarrow} \def\inja{\hookrightarrow}

\newcommand{\lra}{\longrightarrow}
\newcommand{\llra}{\longleftarrow}
\newcommand{\lraf}[1]{\stackrel{#1}{\lra}}
\newcommand{\llaf}[1]{\stackrel{#1}{\llra}}
\newcommand{\ra}{\rightarrow}
\newcommand{\dk}{{\rm dim_{_{k}}}}

\newcommand{\holim}{{\rm Holim}}
\newcommand{\hocolim}{{\rm Hocolim}}
\newcommand{\colim}{{\rm colim\, }}
\newcommand{\limt}{{\rm lim\, }}
\newcommand{\Add}{{\rm Add }}
\newcommand{\Prod}{{\rm Prod }}
\newcommand{\Tor}{{\rm Tor}}
\newcommand{\Cogen}{{\rm Cogen}}
\newcommand{\Tria}{{\rm Tria}}
\newcommand{\Loc}{{\rm Loc}}
\newcommand{\Coloc}{{\rm Coloc}}
\newcommand{\tria}{{\rm tria}}
\newcommand{\Con}{{\rm Con}}
\newcommand{\Thick}{{\rm Thick}}
\newcommand{\thick}{{\rm thick}}
\newcommand{\Sum}{{\rm Sum}}

\begin{center}
{\Large {\bf Mirror-reflective algebras and Tachikawa's second conjecture}}\\ \medskip In memory of Professor Hiroyuki Tachikawa (1930-2022)
\end{center}

\centerline{\textbf{Hongxing Chen}, \textbf{Ming Fang} and \textbf{Changchang Xi}$^*$}

\renewcommand{\thefootnote}{\alph{footnote}}
\setcounter{footnote}{-1} \footnote{ $^*$ Corresponding author.
Email: xicc@cnu.edu.cn; Fax: 0086 10 68903637.}
\renewcommand{\thefootnote}{\alph{footnote}}
\setcounter{footnote}{-1} 
\renewcommand{\thefootnote}{\alph{footnote}}
\setcounter{footnote}{-1} 

\begin{abstract}
Given an algebra  with an idempotent, we introduce two procedures to construct families of new algebras, termed  mirror-reflective algebras and reduced mirror-reflective algebras. We then establish connections among these algebras by recollements of derived module categories. In case of given algebras being gendo-symmetric, we show that the (reduced) mirror-reflective algebras are symmetric and provide new methods to construct systematically both higher dimensional (minimal) Auslander-Gorenstein algebras and gendo-symmetric algebras of higher dominant dimensions. This leads to a new formulation of Tachikawa's second conjecture for symmetric algebras in terms of idempotent stratifications.

\end{abstract}

{\footnotesize\tableofcontents\label{contents}}

\section{Introduction}\label{Introduction}
In the representation theory of algebras, the long-standing and not yet solved Nakayama conjecture says that a finite-dimensional algebra with infinite dominant dimension is self-injective \cite{Nakayama}. This conjecture is related to the so-called Tachikawa's second conjecture \cite{Tac}:

\smallskip
{\bf (TC2)}\;  Let $\Lambda$ be a finite-dimensional self-injective algebra and $M$ a finitely generated $\Lambda$-module. Then $M$ is projective if it is orthogonal, that is,  $\Ext_\Lambda^n(M, M)=0$ for all $n\geq 1$.
\smallskip

By M\"{u}ller's characterization of dominant dimension in \cite{Muller}, {\rm (TC2)} holds for a self-injective algebra $\Lambda$ if and only if the Nakayama conjecture holds for all endomorphism algebras $\End_{\Lambda}(\Lambda\oplus M)$ of finitely generated generators over $\Lambda$. This suggests to consider the algebras $A$ of the form $\End_{\Lambda}(\Lambda\oplus M)$ with $\Lambda$ a self-injective algebra and $M$ an arbitrary $\Lambda$-module. Such algebras are called \emph{Morita algebras} \cite{KY}. In case, $\Lambda$ is symmetric, they are called \emph{gendo-symmetric algebras} \cite{FK16}. In \cite{xc3}, orthogonal generators over a self-injective Artin algebra have been discussed systematically from the viewpoint of recollements of (relative) stable module categories. In particular, it is shown that the Nakayama conjecture holds true for Gorenstein-Morita algebras \cite[Corollary 1.4]{xc3}.

In the present paper, we mainly focus on developing general methods to construct gendo-symmetric algebras. First, for an arbitrary algebra with an arbitrary idempotent, we introduce two inductive procedures to construct families of new algebras, called mirror-reflective algebras and reduced mirror-reflective algebras. In the case of gendo-symmetric algebras, the procedures produce systematically both higher dimensional (minimal) Auslander-Gorenstein algebras and gendo-symmetric algebras of higher dominant dimensions. Second, we show that these families of mirror-reflective algebras are connected by recollements of their derived module categories. Moreover, these recollements are standard in the sense that they are induced from strong idempotent ideals (also called stratifying ideals) of algebras  (see \cite{CPS, CPS2, APT} for definitions). Finally, we give a new formulation of Tachikawa's second conjecture for symmetric algebras in terms of idempotent stratifications.

To state our results more precisely, we first introduce a few terminologies.

Let $A$ be an associative algebra over a commutative ring $k$, $e$ an idempotent of $A$, and $\Lambda:=eAe$. For $\lambda\in Z(\Lambda)$, the center of the algebra $\Lambda$, we introduce an associative algebra  $R(A,e,\lambda)$, called the \emph{mirror-reflective algebra} of $A$ at level $(e,\lambda)$, which has the underlying $k$-module $A \oplus Ae\otimes_\Lambda eA$, such that $Ae\otimes_\Lambda eA$ is an ideal in $R(A,e,\lambda)$ (see Section \ref{sect3.1} for details). The terminology ``mirror-reflective" can be justified by Example \ref{Exmaple1} in Section \ref{MQA}.
Moreover, the $k$-submodule of $R(A,e,\lambda)$
$$ S(A,e,\lambda):=(1-e)A(1-e)\oplus Ae\otimes_\Lambda eA $$
is also an associative algebra with the multiplication induced from the one of $R(A,e,\lambda)$. This algebra is called the \emph{reduced mirror-reflective algebra} of $A$ at level $(e, \lambda)$. It has less simple modules than $R(A,e,\lambda)$ does.
The specialization of  $R(A,e,\lambda)$ and $S(A,e,\lambda)$ at $\lambda=e$ are called the mirror-reflective algebra and reduced mirror-reflective algebra of $A$ at $e$, denoted as $R(A,e)$ and $S(A,e)$, respectively. Moreover, $S(A,e)=e_0R(A,e)e_0$ for an idempotent $e_0$ in $R(A,e)$.

Clearly, each $A$-module is an $R(A,e)$-module via the canonical surjective homomorphism $R(A,e)\to A$ of algebras. Conversely, each $R(A,e)$-module restricts to $A$-module via the canonical inclusion from $A$ into $R(A,e)$. Remark that each module over $(1-e)A(1-e)$ can also be regarded as a module over $S(A,e)$. So we can define two endomorphism algebras associated with $(A,e)$:
$$\mathcal{A}(A,e):=\End_{R(A,e)}\big(R(A,e)\oplus A(1-e)\big),\;\;
\mathcal{B}(A,e):=\End_{S(A,e)}\big(S(A,e)\oplus (1-e)A(1-e)\big).$$

Now, assume that all algebras considered are finite-dimensional associative $k$-algebras with identity over a field $k$.
Further, assume that $A$ is a gendo-symmetric algebra and $e$ is an idempotent of $A$ such that $Ae$ is a faithful, projective-injective $A$-module. In this case, we write $(A, e)$ for the gendo-symmetric algebra $A$. If $e'$ is another idempotent of $A$ such that $Ae'$ is a faithful, projective-injective $A$-module, then $R(A,e)\simeq R(A,e')$ as algebras (see Lemma \ref{generator}(1)). Hence, up to isomorphism of algebras, we can write $R(A)$ for $R(A,e)$ without referring to $e$, and call it the \emph{mirror-reflective algebra} of the gendo-symmetric algebra $A$.

Our first result reveals homological properties of mirror-reflective  algebras of gendo-symmetric algebras.
Recall that $A$ is called an \emph{$n$-Auslander algebra} ($n\ge 0$) if $\gd(A)\leq n+1\leq\dm(A)$; an \emph{$n$-minimal Auslander-Gorenstein algebra} if $\id({_A}A)\leq n+1\leq\dm(A)$ (see \cite{Au71,Iyama07,IS,CK}), where $\gd(A)$, $\dm(A)$ and $\id(_AA)$ denote the global, dominant and left injective dimensions of an algebra $A$, respectively. Clearly, $n$-Auslander algebras are exactly $n$-minimal Auslander-Gorenstein algebras of finite global dimension (see Subsection \ref{FDDD}).

\begin{Theo} \label{Main properties}
Let $(A,e)$ be a gendo-symmetric algebra.  Then

$(1)$ $R(A,e,\lambda)$ is a symmetric algebra for $\lambda$ in the center of $eAe$.

$(2)$ $\min\{\dm(\mathcal{A}(A,e)), \dm(\mathcal{B}(A,e))\}$
$\geq \dm(A)+2.$

$(3)$ $\mathcal{A}(A,e)$ is $(2n+3)$-minimal Auslander-Gorenstein if $A$ is $n$-minimal Auslander-Gorenstein with $n>0$. Further, $\mathcal{A}(A,e)$ is an $(2n+3)$-Auslander algebra if $A$ is an $n$-Auslander algebra.
\end{Theo}

According to Theorem \ref{Main properties}, $\mathcal{A}(A,e)$ and $\mathcal{B}(A,e)$ are gendo-symmetric. Thus the construction of (reduced) mirror-reflective algebras can be done iteratively. Starting with a fixed gendo-symmetric algebra $(A,e)$, we may define $4$ families of algebras: $R_n$, $S_n$, $A_n$ and $B_n$ for $n>0$ (see Section \ref{IMSAT} for details). They are called the $n$-th mirror reflective, $n$-th reduced mirror-reflective, $n$-th gendo-symmetric and $n$-th reduced gendo-symmetric algebras of $(A, e)$, respectively. By Theorem \ref{Main properties}(2), these algebras have higher homological dimension: $\dm(A_{n+1})\geq \dm(A_n)+2$ and $\dm(B_{n+1})\geq \dm(B_n)+2$. Thus
$\min\{\dm(A_n),\dm(B_n)\}\geq \dm(A)+2(n-1)\geq 2n.$ Moreover, they are connected by derived recollements, as is shown in our second result below. Here, $\Df{A}$ and $\D{A}$ denote the bounded-above and unbounded derived categories of $A$, respectively.

\begin{Theo} \label{Recollement} Let $(A,e)$ be a gendo-symmetric algebra and $n$ a positive integer. Then the following hold.

$(1)$ There exist recollements of bounded-above derived categories of algebras:

$$\xymatrix{\Df{A_n}\ar[r]&\Df{A_{n+1}}\ar[r]\ar@/^1.2pc/[l]\ar@/_1.2pc/[l]
&\Df{A_n} \quad and\ar@/^1.2pc/[l]\ar@/_1.2pc/[l]}\vspace{0.2cm}\quad
\xymatrix{\Df{B_0}\ar[r]&\Df{B_{n+1}}\ar[r]\ar@/^1.2pc/[l]\ar@/_1.2pc/[l]
&\Df{B_n}\ar@/^1.2pc/[l]\ar@/_1.2pc/[l]}\vspace{0.2cm}$$
with $B_0:=(1-e)A(1-e)$.

$(2)$ Let $R_0=S_0:=eAe$. If $\dm(A)=\infty$, then there are recollements of unbounded derived categories of algebras
induced by strong idempotent ideals:

$$\xymatrix{\D{A_n}\ar[r]&\D{R_n}\ar[r]\ar@/^1.2pc/[l]\ar@/_1.2pc/[l]
&\D{R_{n-1}} \quad and \ar@/^1.2pc/[l]\ar@/_1.2pc/[l]}\vspace{0.2cm}\quad \xymatrix{\D{B_0}\ar[r]&\D{S_n}\ar[r]\ar@/^1.2pc/[l]\ar@/_1.2pc/[l]
&\D{S_{n-1}}.\ar@/^1.2pc/[l]\ar@/_1.2pc/[l]}\vspace{0.2cm}$$
\end{Theo}

Motivated by Theorem \ref{Recollement}(2), we introduce the \emph{stratified dimension} of an algebra. It measures how many steps an algebra can be stratified by its nontrivial strong idempotents (see Definition \ref{Sdim}), or equivalently, the derived category of the algebra can be stratified by nontrivial standard recollements of derived module categories. We also define the \emph{stratified ratio} of an algebra to be the ratio of its stratified dimension to the number of isomorphism classes of simple modules (see Definition \ref{Sr}).

Our third result establishes a connection between ${\rm (TC2)}$ and stratified dimensions of algebras.

\begin{Theo}  \label{Tachikawa's second conjecture}
The following are equivalent for a field $k$.

$(1)$ Tachikawa's second conjecture holds for all  symmetric algebras over $k$.

$(2)$ Each indecomposable symmetric algebra over $k$ has no nontrivial strong idempotent ideals.

$(3)$ The supreme of stratified ratios of all indecomposable symmetric algebras over $k$ is less than $1$.
\end{Theo}

As a consequence of Theorem \ref{Tachikawa's second conjecture},  we shall provide a sufficient condition for  {\rm (TC2)} on symmetric algebras to hold true. An algebra $\Lambda$ is said to be \emph{derived simple} if its unbounded derived category $\D{\Lambda}$ admits no nontrivial recollements of unbounded derived module categories of algebras. Examples of derived simple algebras include local algebras, blocks of groups algebras and some indecomposable algebras with two simple modules. One should not confuse the notion of derived simple algebras with the one of \emph{$\mathscr{D}^{\rm b}(\rm mod)$-derived simple} algebras in the sense that the bounded derived categories (of finitely generated modules) are not  nontrivial recollements of bounded derived categories of any algebras (see \cite{LY}). Derived simple algebras are  $\mathscr{D}^{\rm b}(\rm mod)$-derived simple, but the converse is not true in general. By \cite[Theorem 3.2]{LY}, each indecomposable symmetric algebra is always $\mathscr{D}^{\rm b}(\rm mod)$-derived simple.

If an algebra $\Gamma$ has a nontrivial strong idempotent ideal generated by an idempotent $f$ (see Definition \ref{Strong} below), then there is a nontrivial recollement $(\D{\Gamma/\Gamma f\Gamma},\D{\Gamma}, \D{f\Gamma f})$. Thus we obtain the following corollary immediately from Theorem \ref{Tachikawa's second conjecture}.

\begin{Koro}\label{Derived simple}
If each indecomposable symmetric algebra over a field $k$ is derived simple, then Tachikawa's second conjecture holds for all symmetric algebras over $k$.
\end{Koro}

The paper is structured as follows. In Section $2$, we recall the definitions of strong idempotents, recollements as well as higher Auslander-Gorenstein and Auslander algebras. Also, we present the definitions of  stratified dimensions and ratios of algebras (see Definitions \ref{Sdim} and \ref{Sr}), respectively. In Section $3$, we define (reduced) mirror-reflective algebras by reflecting a left (or right) ideal generated by an idempotent element. We then study the derived module categories of these algebras. Also, an explicit description of mirror-reflective algebras is demonstrated by quivers with relations. This explains visually the terminology of mirror-reflective algebras. In Section $4$, we first show Theorems \ref{Main properties} and \ref{Recollement}. This relies on the fact that mirror-reflective algebras of gendo-symmetric algebras at any levels are symmetric (see Proposition \ref{Mirror-symmetric}). By iteration of forming (reduced) mirror-reflective algebras from a gendo-symmetric algebra, a series of recollements of derived module categories is established. This not only gives a proof of Theorem \ref{Tachikawa's second conjecture}, but also shows the relation between the numbers of simple modules over different mirror-reflective algebras (see Corollary \ref{FKK}(2)-(3)).

\section{Preliminaries}
Throughout the paper, $k$ denotes a commutative ring, and all algebras considered are associative $k$-algebras with identity.

Let $A$ be a $k$-algebra. We denote by $A\Modcat$ the category of all left $A$-modules, and by $A\modcat$ the full subcategory of $A\Modcat$  consisting of finitely generated $A$-modules. The \emph{global dimension} of $A$, denoted by $\gd(A)$, is defined to be the supreme of projective dimensions of all $A$-modules. The \emph{finitistic dimension} of $A$, denoted by $\fd(A)$, is defined to be the supreme of projective dimensions of those $A$-modules which have a finite projective resolution by finitely generated projective modules.  The \emph{ projective and injective dimensions} of an $A$-module $M$ are denoted by $\pd(_AM)$ and $\id({_A}M)$, respectively. If $f:X\to Y$ and $g:Y\to Z$ are homomorphisms of $A$-modules, we write $fg$ for the composite of $f$ with $g$, and $(x)f$ for the image of $x\in X$ under $f$.

For an additive category $\mathcal{C}$, let $\C{\mathcal{C}}$ denote the category of all complexes over $\mathcal{C}$ with chain maps, and $\K{\mathcal{C}}$ the homotopy category of $\C{\mathcal{C}}$. We denote by $\Cb{\mathcal{C}}$ and $\Kb{\mathcal{C}}$ the full subcategories of $\C{\mathcal{C}}$ and
$\K{\mathcal{C}}$ consisting of bounded complexes over $\mathcal{C}$, respectively. When $\mathcal{C}$ is abelian, the \emph{(unbounded) derived category} of $\mathcal{C}$ is denoted by $\D{\mathcal{C}}$, which is the localization of $\K{\mathcal C}$ at all quasi-isomorphisms. The full subcategory of $\D{\mathcal{C}}$
consisting of bounded-above complexes over $\mathcal{C}$ is denoted by $\Df{\mathcal{C}}$. As usual, we simply write
$\K{A}$ for $\K{A\Modcat}$, $\D{A}$ for $\D{A\Modcat}$, and $\Df{A}$ for $\Df{A\Modcat}$. Also, we identify $A\Modcat$ with the
full subcategory of $\D{A}$ consisting of all stalk complexes in degree zero.

\subsection{Standard recollements and stratified dimensions}\label{SRII}
In this section, we start with recalling recollements of triangulated categories, introduced by Beilinson, Bernstein and Deligne in \cite{BBD}, and introduce the notion of stratified dimensions of algebras.

\begin{Def}\label{def01} \rm
Let  $\mathcal{T}$, $\mathcal{T'}$ and $\mathcal{T''}$ be
triangulated categories. $\mathcal{T}$ admits a
\emph{recollement} of $\mathcal{T'}$ and $\mathcal{T''}$ (or there is a recollement among $\mathcal{T''}, \mathcal{T}$ and $\mathcal{T}'$) if there
are six triangle functors
$$\xymatrix@C=1.5cm{\mathcal{T''}\ar^-{i_*=i_!}[r]&\mathcal{T}\ar^-{j^!=j^*}[r]
\ar^-{i^!}@/^1.2pc/[l]\ar_-{i^*}@/_1.6pc/[l]
&\mathcal{T'}\ar^-{j_*}@/^1.2pc/[l]\ar_-{j_!}@/_1.6pc/[l]}$$ among the three categories such
that the $4$ conditions are satisfied:

$(1)$ $(i^*,i_*),(i_!,i^!),(j_!,j^!)$ and $(j^*,j_*)$ are adjoint
pairs.

$(2)$ $i_*,j_*$ and $j_!$ are fully faithful functors.

$(3)$ $j^! i_!=0$ (and thus also $i^!j_*=0$ and $i^*j_!=0$).

$(4)$ For an object $X\in\mathcal{T}$, there are triangles
$i_!i^!(X)\to  X\to j_*j^*(X)\to i_!i^!(X)[1]$ and $ j_!j^!(X)\to
X\to i_*i^*(X)\to j_!j^!(X)[1]$ induced by the adjunctions of counits and units,
where [$1$] is the shift functor of $\mathcal{T}$.
\end{Def}

Quasi-hereditary algebras, introduced by Cline, Parshall and Scott (see \cite{CPS,CPS2}), provide a special class of recollements of derived module categories. Recall that a\emph{ heredity ideal} of an algebra $A$ is an ideal $I$ such that $(i)$ $I$ is idempotent (i.e. $I^2=I$), $(ii)$ $_AI$ is projective as an $A$-module and $(iii)$ $\End_A(_AI)$ is semisimple. For such an ideal $I$, there holds always $\Ext_{A/I}^i(X,Y) \simeq \Ext_A^i(X,Y)$ for all modules $X,Y$ over $A/I$ and for all $i\ge 0$. A slight generalisation of heredity ideals is the $n$-idempotent ideals defined in \cite{APT}.

\begin{Def}{\rm \cite{APT}}\label{Strong}
Let $A$ be an algebra, $I$ an ideal of $A$, and $n$ a positive integer. The ideal $I$ of $A$ is said to be \emph{$n$-idempotent} if, for  $X,Y\in (A/I)\Modcat$, the canonical homomorphism $\Ext_{A/I}^i(X,Y) \to \Ext_A^i(X,Y)$ of $k$-modules is an isomorphism for all $1\leq i\leq n$.

The ideal $I$ is said to be \emph{strong idempotent} if $I$ is $n$-idempotent for all $n\geq 1$.
In this case, if $I=AeA$ for an idempotent $e\in A$, then $e$ is called a \emph{strong idempotent} of $A$.
\end{Def}

By a trivial strong idempotent of $A$ we mean the idempotent $0$ or an idempotent $e$ with $AeA=A$.
Note that an ideal $I$ is $1$-idempotent if and only if $I$ is idempotent. Moreover, strong idempotent ideals are closely related to homological ring epimorphisms. Recall that a ring homomorphism $\lambda:A\to B$ is called a \emph{homological ring epimorphism} if the multiplication map $B\otimes_AB\to B$ is an isomorphism and $\Tor_i^A(B,B)=0$ for all $i\geq 1$. This is equivalent to saying that the derived restriction functor $D(\lambda_*):\D{B}\to\D{A}$, induced by the restriction functor $\lambda_*: B\Modcat\to A\Modcat$, is fully faithful.
Note that strong idempotent ideals $I$ are also called homological ideals, that is, the canonical surjection $A\to A/I$ is a homological ring epimorphism.

Let us emphasize that strong idempotent ideals generated by idempotents are exactly \emph{stratifying ideals} introduced in \cite[Definition 2.1.1]{CPS2}.

\begin{Lem}{\rm \cite{APT}}\label{strong idempotent}
Let $I=AeA$ for an idempotent $e$ in $A$.

$(1)$ Let $n$ be a positive integer. Then $I$ is $(n+1)$-idempotent if and only if
the multiplication map $$Ae\otimes_{eAe}eA\lra I, \;\;ae\otimes eb\mapsto aeb,\;\, a,b\in A$$ is an isomorphism of $A$-$A$-bimodules and $\Tor_i^{eAe}(Ae,eA)=0$ for all $1\leq i\leq n-1$.

$(2)$ If $I$ is $2$-idempotent, then
$$\sup\{n\in\mathbb{N}\mid \Ext_A^i(A/I, A/I)=0, 1\leq i\leq n\}\geq \sup\{n\in\mathbb{N}\mid\Tor_i^{eAe}(Ae,eA)=0, 1\leq i\leq n\}+2.$$
\end{Lem}

{\it Proof.} $(1)$ Although all the results in \cite{APT} are stated for finitely generated modules over Artin algebras, many of them such as Theorem 2.1, Lemma 3.1 and Propositions 2.4 and 3.7(b) hold for arbitrary modules over rings if we modify $\mathbf{P}_n$ in \cite[Definition 2.3]{APT}
as follows:

Let $\mathbf{P}_n(Ae)$ be the full subcategory of $A\Modcat$ consisting of all modules $X$ such that there is an exact
sequence $P_n\to\cdots\to P_1\to P_0\to X\to 0$ of $A$-modules with $P_i\in\Add(Ae)$ for $0\leq i\leq n$, where $\Add(Ae)$ is the full subcategory of $A\Modcat$ consisting of direct summands of direct sums of copies of $Ae$.

By \cite[Theorem 2.1]{APT}, $I:=AeA$ is $(n+1)$-idempotent if and only if $I\in\mathbf{P}_n(Ae)$.
In particular, $I$ is $2$-idempotent if and only if $I\in\mathbf{P}_1(Ae)$. By \cite[Lemma 3.1]{APT},
the adjoint pair $(Ae\otimes_{eAe}-, \Hom_A(Ae,-))$ between $(eAe)\Modcat$ and $A\Modcat$ induces additive equivalences between
$(eAe)\Modcat$ and $\mathbf{P}_1(Ae)$. Note that $\Hom_A(Ae, I)\simeq eI=eA$.
Thus $I\in\mathbf{P}_1(Ae)$ if and only if the multiplication map $Ae\otimes_{eAe}eA\to AeA$
is an isomorphism of $A$-$A$-bimodules. Now, assume that $I$ is $2$-idempotent. By \cite[Proposition 3.7(b)]{APT}, $I\in\mathbf{P}_n(Ae)$
if and only if $\Tor_i^{eAe}(Ae, eA)=0$ for all $1\leq i< n$. This shows $(1)$.

$(2)$ If $I$ is $(n+1)$-idempotent, then $\Ext_A^i(A/I, A/I)\simeq\Ext_{A/I}^i(A/I, A/I)=0$ for all $1\leq i\leq n+1$. Now, $(2)$ follows from $(1)$. $\square$

\begin{Koro}\label{HID}
$(1)$ Let $e$ and $f$ be idempotents of $A$ such that $ef=e=fe$. If $AeA$ is an $(n+1)$-idempotent ideal of $A$ for a positive integer $n$, then $fAeAf$ is
an $(n+1)$-idempotent ideal of $fAf$. In particular, if $e$ is a strong idempotent of $A$, then it is also a strong idempotent of $fAf$.

$(2)$ Let $\{e,e_1,e_2\}$ be a set of orthogonal idempotents of $A$ such that $e$ is a strong idempotent of $A$. Define $f:=e+e_1$, $g:=e+e_1+e_2$ and $\overline{A}:=A/AeA$. Let $\overline{f}:=f+AeA$ denote the image of $f$ in $\overline{A}$. If $\overline{f}$ is a strong idempotent of $\overline{g}\overline{A}\overline{g}$, then $f$ is a strong idempotent of $gAg$.
\end{Koro}

{\it Proof.} $(1)$ Transparently, $e\in fAf$, $efAfe=eAe$, $fAeAfe=fAe$ and $efAeAf=eAf$. If $Ae\otimes_{eAe}eA\simeq AeA$, then $fAe\otimes_{eAe}eAf\simeq fAeAf$.
Since $Ae=fAe\oplus (1-f)Ae$ and $eA=eAf\oplus eA(1-f)$, we see that the abelian group $\Tor_i^{eAe}(fAe, eAf)$ is a direct summand of $\Tor_i^{eAe}(Ae, eA)$
for $i\in \mathbb{N}$. Now, (1) follows from Lemma \ref{strong idempotent}(1).

$(2)$ Clearly, $AeA\subseteq AfA\subseteq AgA$, and $\overline{g}\overline{A}\overline{g}\simeq gAg/gAeAg$ and $\overline{g}\overline{A}\overline{g}/\overline{g}\overline{A}\overline{f}\overline{A}\overline{g}\simeq gAg/gAfAg$ as algebras.
Suppose that $\overline{f}$ is a strong idempotent of $\overline{g}\overline{A}\overline{g}$. Then the canonical surjection $\pi_2: gAg/gAeAg\to gAg/gAfAg$ is homological. Since $e$ is a strong idempotent of $A$ and $ge=e=eg$, the canonical surjection $\pi_1: gAg\to gAg/gAeAg$ is also homological by (1). Observe that compositions of homological ring epimorphisms are again homological ring epimorphisms. Thus $\pi_1\pi_2: gAg\to gAg/gAfAg$ is homological. This implies that $f$ is a strong idempotent in $gAg$. $\square$

\medskip

Let $e=e^2\in A$. If $AeA$ is a strong idempotent ideal in $A$, then the recollement of derived module categories of algebras:
$$\xymatrix{\mathscr{D}(A/AeA)\ar[r]&\mathscr{D}(A)\ar[r]\ar@/^1.2pc/[l]\ar@/_1.2pc/[l]
&\mathscr{D}(eAe)\ar@/^1.2pc/[l]\ar@/_1.2pc/[l]}\vspace{0.2cm}.$$

\noindent is called a \emph{standard recollement} induced by $AeA$.
If ${_A}AeA$ or $AeA_A$ is projective (for example, $AeA$ is a heredity ideal in $A$), then the ideal $AeA$ is strong idempotent. In the case that ${_A}AeA$ is projective, the recollement restricts to a recollement $(\Df{A/AeA}, \Df{A}, \Df{eAe})$ of bounded-above derived categories.

A general method is given for constructing finitely generated (one-sided) projective idempotent ideals of the endomorphism algebras of objects in additive categories (see \cite[Lemmas 3.2 and 3.4]{xc4}). This implies the following.

\begin{Lem}\label{Projective}
Suppose that $R$ is an algebra and $I$ is an ideal of $R$.

$(1)$ Let $A:=\End_R(R\oplus \; R/I)$ and $e^2=e\in A$ correspond to the direct summand $R/I$ of the $R$-module
$R\oplus \; R/I$. Then $AeA{_A}$ is finitely generated and projective, and there is a recollement $(\D{R/{\rm Ann}_{R\opp}(I)},$ $ \D{A}, \D{R/I})$, with ${\rm Ann}_{R\opp}(I):=\{r\in R\mid Ir=0\}$.

$(2)$ Let $B:=\End_R(R\oplus I)$ and $f=f^2\in B$ correspond to the direct summand $I$ of the $R$-module $R\oplus I$.  If $I$ is idempotent, then ${_B}BfB$ is finitely generated and projective, and there is a recollement $(\D{R/I}, \D{B}, \D{\End_R(I)})$.
\end{Lem}

Another way to produce finitely generated projective ideals comes from Morita context algebras, as explained below.

Let $R$ be an algebra and let $I$ and $J$ be ideals of $R$ with $IJ=0$. Define $$M_l(R,I,J):=\left(\begin{array}{lc} R & I\\
R/J & R/J\end{array}\right)\quad(\;\mbox{respectively,} \;\;M_r(R,I,J):=\left(\begin{array}{lc} R & R/I\\
J & R/I\end{array}\right)\;)$$ which is the \emph{Morita context algebra} with the bimodule homomorphisms given by the canonical ones:
$$
I\otimes_{R/J}(R/J)\simeq I\hookrightarrow R,\quad (R/J)\otimes_R I\simeq I/JI\twoheadrightarrow (I+J)/J\hookrightarrow R/J
$$
$
(\;\mbox{respectively,}\;\;(R/I)\otimes_{R/I}J\simeq J\hookrightarrow R,\quad J\otimes_R(R/I)\simeq J/JI\twoheadrightarrow (I+J)/I\hookrightarrow R/I).
$
Note that $M_r(R,I, {\rm Ann}_{R\opp}(I))\simeq \End_R(R\oplus R/I)$ as algebra. Moreover, if ${_R}R$ is injective and $I^2=I$, then
$M_l(R,I,{\rm Ann}_{R\opp}(I))\simeq \End_R(R\oplus I)$ as algebras. This is due to $\Hom_R(I,R/I)=0.$

Let $$e:=\left(\begin{array}{lc} 0 & 0\\
0& 1+J\end{array}\right)\in M_l(R,I,J),\;\; f:=\left(\begin{array}{lc} 0 & 0\\
0 & 1+I\end{array}\right)\in M_r(R,I,J).$$
Then the next lemma is easy to verify.
\begin{Lem}\label{LRMOR}
Let $A:=M_l(R,I,J)$ and $B:= M_r(R,I,J)$. Then ${_A}AeA$ and $BfB_B$ are finitely generated and projective.
Moreover, there are recollements $(\D{R/I}, \D{A}, \D{R/J})$ and $(\D{R/J}, \D{B}, \D{R/I})$.
\end{Lem}

Now, we introduce stratified dimensions of algebras, which measure how many steps the given algebras can be stratified by their nontrivial strong idempotents.

\begin{Def}\label{Sdim}
By an \emph{idempotent stratification} of length $n$ of an algebra $A$, we mean a set $\{e_i\mid 0\leq i\leq n\}$ of nonzero (not necessarily primitive) orthogonal idempotents of $A$ satisfying the conditions:

$(a)$ $1=\sum_{j=0}^{n}e_j$ and $e_{i+1}\notin Ae_{\leq i}A$ (or equivalently, $Ae_{\leq i}A \subsetneq Ae_{\leq (i+1)}A$) for all $0\leq i\leq n-1$, where $e_{\leq m}:=\sum_{j=0}^{m}e_j$ for each $0\leq m\leq n$; and

$(b)$ $e_{\leq i}$ is a strong idempotent of the algebra $e_{\leq (i+1)}Ae_{\leq (i+1)}$ for each $0\leq i\leq n-1$.

\smallskip
The \emph{stratified dimension} of $A$, denoted by $\sd(A)$, is defined to be the supreme of the lengths of all idempotent stratifications of $A$.
\end{Def}

Clearly, $\sd(A)=0$ if and only if $A$ has no nontrivial strong idempotent ideals. If $\sd(A)=n>0$, then
there are iterated nontrivial standard recollements
$\big(\D{A_{i}/I_i},\D{A_i}, \D{A_{i-1}}\big)$ for all $1\leq i<n+1$, where
$A_0:=e_0Ae_0$, $A_i:=e_{\leq i}Ae_{\leq i}$ and $I_i:=e_{\leq i}Ae_{\leq (i-1)}Ae_{\leq i}$ in Definition \ref{Sdim}.
Moreover, for any two algebras $A_1$ and $A_2$, $\sd(A_1\times A_2)=\sd(A_1)+\sd(A_2)+1.$
This implies that the stratified dimension of the direct product of $\mathbb{N}$-copies of a field $k$ is infinite.

Stratifications of algebras in the sense of Cline, Parshall and Scott are idempotent stratifications. But the converse is not true.
Following \cite[Chapter 2]{CPS2}, a \emph{stratification} of length $(n+1)$ of an algebra $A$ is a chain of ideals,
$0=U_{-1}\subsetneq U_0\subsetneq U_1\subsetneq\cdots\subsetneq U_{n-1}\subsetneq U_n=A$, generated by idempotents such that $U_i/U_{i-1}$ is a strong idempotent ideal of $A/U_{i-1}$ for  $0\leq i\leq n$. In this case, $A$ is said to be \emph{CPS-stratified}.
If $\{e_i\mid 0\leq i\leq n\}$ is a complete set of nonzero primitive orthogonal idempotents of $A$ and $U_i=Ae_{\leq i}A$ for  $0\leq i\leq n$, then $A$ is called a \emph{fully CPS-stratified} algebra. Standardly stratified algebras with respect to an order of simple modules are fully CPS-stratified.

\begin{Lem}\label{homological}
Let $\{e_i\mid 0\leq i\leq n\}$ be a set of nonzero orthogonal idempotents of $A$ satisfying the condition $(a)$ in Definition \ref{Sdim}.
Define $U_i:=Ae_{\leq i}A$ for $0\leq i\leq n$ and $U_{-1}:=0$. If $U_i/U_{i-1}$ is a strong idempotent ideal of $A/U_{i-1}$ for $0\leq i\leq n$, then
the condition $(b)$ in Definition \ref{Sdim} holds.
 \end{Lem}

{\it Proof.} Since $U_i/U_{i-1}$ is a strong idempotent ideal of $A/U_{i-1}$ by assumption, the canonical surjection $A/U_{i-1}\to A/U_i$ is homological.
As the composition of homological ring epimorphisms is still a homological ring epimorphism, the canonical surjection $A\to A/U_i$ is  homological.
This implies that $e_{\leq i}$ is a strong idempotent of $A$. By Corollary \ref{HID}, $e_{\leq i}$ is a strong idempotent of $e_{\leq (i+1)}Ae_{\leq (i+1)}$.
Thus Definition \ref{Sdim}(b) holds. $\square$

\smallskip
For an Artin algebra $A$, we denote by $\#(A)$ the number of isomorphism classes of simple $A$-modules.

\begin{Prop}\label{upper}
Let $A$ be an Artin algebra over a commutative Artin ring $k$.  Then

$(1)$ $\sd(A)\leq \#(A)-1$.

$(2)$ If $A$ has a stratification of length $n+1$ with $n\in\mathbb{N}$, then $\sd(A)\geq n$. In particular, if $A$ is a fully CPS-stratified algebra, then $\sd(A)=\#(A)-1$.

$(3)$ If $\sd(A)\geq 1$, then
$\sd(A)=\sup_{e\in A}\{\sd(eAe)+\sd(A/AeA)+1\}$,
where $e$ runs over all nonzero strong idempotents of $A$ with $AeA\neq A$.

$(4)$ If $k$ is a field and $B$ is a finite-dimensional $k$-algebra, then
$$\sd(A\otimes_kB)\geq (\sd(A)+1)(\sd(B)+1)-1.$$
\end{Prop}

{\it Proof.} $(1)$ This is clear by Definition \ref{Sdim}(a).

$(2)$ The first part of $(2)$ follows from Lemma \ref{homological}.  If $A$ is a fully CPS-stratified algebra, then
it has a stratification of length $\#(A)-1$. By $(1)$, we obtain $\sd(A)=\#(A)-1$.

$(3)$ An Artin algebra always has only finitely many nonisomorphic, indecomposable, finitely generated projective modules. This implies

$(\ast)$ If $f$ is an idempotent of $A$ and $I$ is an idempotent ideal of $A$ such that $AfA\subseteq I$, then there is an idempotent $f'$ of $A$ which is orthogonal to $f$ such that $I=A(f+f')A$.

Now, let $n:=\sd(A)\geq 1$. On the one hand, since $e_{\leq n-1}$ in Definition \ref{Sdim}(b) is a strong idempotent of $A$, we have $\sd(A)=\sd(e_{\leq n-1}Ae_{\leq n-1})+1$ and $\sd(A/Ae_{\leq n-1}A)=0$ by $(\ast)$ and Corollary \ref{HID}(2).  On the other hand, for each nonzero strong idempotent $e$ of $A$ with $AeA\neq A$, it follows again from $(\ast)$ and Corollary \ref{HID}(2) that $\sd(eAe)+\sd(A/AeA)+1\leq n$. Thus $(3)$ holds.

$(4)$ Let $m:=\sd(B)$ and $\ell:=n+m$. If $\ell=0$ (i.e. $n=0=m$), then the inequality holds obviously.  Let $\ell\geq 1$.  Without loss of generality, suppose $n\geq 1$. By the proof of $(3)$, there is a nonzero strong idempotent $e$ of $A$  with $AeA\neq A$ such that $\sd(eAe)=n-1$ and $\sd(A/AeA)=0$. Then  the canonical surjection $\pi: A\to A/AeA$ is homological.  Note that, for homological epimorphisms $\lambda_i: R_i\to S_i$ of  algebras over the field $k$ with $i=1,2$,  the algebra homomorphism $\lambda_1\otimes_k\lambda_2: R_1\otimes_k R_2\to S_1\otimes_k S_2$ is still a homological ring epimorphism.  This is due to the isomorphism
$$\Tor_j^{R_1\otimes_k R_2}(S_1\otimes_k S_2, S_1\otimes_k S_2)\simeq\bigoplus_{p+q=j}\Tor_p^{R_1}(S_1,S_1)\otimes_k\Tor^{R_2}_q(S_2, S_2) \; \mbox{ for all } j\in\mathbb{N}.$$
Now, let $C:=A\otimes_kB$ and $e':=e\otimes 1\in C$. Then the surjection $\pi\otimes 1: C\to A/AeA\otimes_k B$  is homological.  Clearly,  there are algebra isomorphisms $(A/AeA)\otimes_k B \simeq C/(AeA\otimes_kB)\simeq C/Ce'C$. It follows that the canonical surjection $C\to C/Ce'C$ is homological, and therefore $e'$ is a nonzero strong idempotent of $C$ with $Ce'C\neq C$. By $(3)$,  $\sd(C)\geq \sd(eAe\otimes_k B)+\sd((A/AeA)\otimes_k B)+1$. Moreover, by induction, $ \sd(eAe\otimes_k B)\geq (\sd(eAe)+1)(\sd(B)+1)-1$ and $\sd((A/AeA)\otimes_k B)\geq \sd(B)$.  Thus  $\sd(C)\geq (n+1)(m+1)-1$.  $\square$

\begin{Def}\label{Sr}
Let $A$ be an Artin algebra over a commutative Artin ring $k$.  The rational number  $\frac{\sd(A)}{\#(A)}$ is called
the \emph{stratified ratio} of $A$ and denoted by $\sr(A)$.
\end{Def}

By Proposition \ref{upper}(1), $\sr(A)\in \mathbb{Q}\cap [0,1)$. Let $A^n$ denote the product of $n$-copies of $A$. Then
$$\lim\limits_{n\to\infty}\sr(A^n)=\lim\limits_{n\to\infty}\frac{n(\sd(A))+n-1}{n\;\#(A)}=\frac{\sd(A)+1}{\#(A)}\leq 1.$$
In particular, if $\sd(A)=\#(A)-1$ (for example, $A$ is quasi-hereditary or local), then $\lim\limits_{n\to\infty}\sr(A^n)=1$. In Section \ref{MSAGSA},  for a gendo-symmetric algebra with infinite dominant dimension, we construct a series of indecomposable symmetric algebras $S_n$ such that $\lim\limits_{n\to\infty}\sr(S_n)=1$ (see Corollary \ref{DMSD}).

\subsection{Dominant dimensions  and gendo-symmetric algebras }\label{FDDD}

Let $A$ be a finite-dimensional algebra over a field $k$.

\begin{Def}
The \emph{dominant dimension} of an algebra $A$, denoted by $\dm(A)$, is the maximal natural number $n$ or $\infty$ such that the first $n$ terms $I_0, I_1, \cdots, I_{n-1}$ in a minimal injective resolution
$0\rightarrow {}_AA\rightarrow I^0\rightarrow I^1\rightarrow \cdots\rightarrow I^{i-1}\rightarrow I^i\rightarrow \cdots$ of $A$
are projective.
\end{Def}

Recall that a module $M\in A\modcat$ is called a \emph{generator} if $A\in\add(M)$; a \emph{cogenerator} if $D(A)\in\add(M)$;  a \emph{generator-cogenerator} if it is a generator and cogenerator. By \cite[Lemma 3]{Muller}, if $_AM\in A\modcat$ is a generator-cogenerator, then $\dm(\End_A(M))=\sup\{n\in\mathbb{N}\mid \Ext_A^i(M,M)=0, 1\leq i\leq n\}+2.$

Algebras of the form $\End_A(A\oplus M)$ with $A$ an algebra and $M$ an $A$-module has double centralizer property and has been studied for a long time. Following  \cite{FK16}, such an algebra is called a \emph{gendo-symmetric} algebra if the algebra $A$ is symmetric. Note that if $A$ is a symmetric algebra, then so is $eAe$ for $e=e^2\in A$.

\begin{Lem}{\rm \cite[Theorem 3.2]{FK11}}\label{GS}
The following are equivalent for an algebra $A$ over a field.

$(1)$ $A$ is a gendo-symmetric algebra.

$(2)$ $\dm(A)\geq 2$ and $D(Ae)\simeq eA$ as $eAe$-$A$-bimodules, where $e\in A$ is an idempotent such that $Ae$  is a faithful projective-injective $A$-module.

$(3)$ $\Hom_A(D(A),A)\simeq A$ as $A$-$A$-bimodules.

$(4)$ $D(A)\otimes_AD(A)\simeq D(A)$ as $A$-$A$-bimodules.
\end{Lem}

In the rest of the paper, we always write $(A, e)$ for a gendo-symmetric
algebra with $e$ an idempotent in $A$ such that $Ae$ is a \emph{faithful projective-injective} $A$-module. Note that $\add(Ae)$ coincides with the full subcategory of $A\modcat$ consisting of projective-injective $A$-modules.

An algebra $A$ is called an \emph{Auslander algebra} if $\gd(A)\leq 2\leq\dm(A)$. This is equivalent to saying that $A$ is the endomorphism algebra of an additive generator of a representation-finite algebra over a field (see \cite{Au71}). A generalization of Auslander algebras is the so-called $n$-Aulslander algebras. Let $n$ be a positive integer. Following \cite{Au71,Iyama07,IS}, $A$ is called an \emph{$n$-Auslander algebra} if $\gd(A)\leq n+1\leq\dm(A)$; an \emph{$n$-minimal Auslander-Gorenstein algebra} if $\id({_A}A)\leq n+1\leq\dm(A)$. Clearly, $n$-Auslander algebras are $n$-minimal Auslander-Gorenstein, while $n$-minimal Auslander-Gorenstein algebras of finite global dimension are $n$-Auslander.  Moreover, these algebras can be characterized in terms of left or right perpendicular categories. For each $M\in A\modcat$ and $m\in\mathbb{N}$, we define
$$^{\perp_{m}}M:=\{X\in A\modcat\mid\Ext_A^i(X,M) = 0, 1\leq i\leq m~\},\;\;M^{\perp_m}:=\{X\in A\modcat\mid\Ext_A^i(M,X) = 0, 1\leq i\leq m~\}.$$
Recall that an $A$-module $N$ is said to be \emph{maximal $(n-1)$-orthogonal} or \emph{$n$-cluster tilting} if $\add({_\Lambda}N)= {}^{\perp_{n-1}}N=N^{\perp_{n-1}}$. A generator-cogenerator $A$-module $M$ is said to be \emph{$(n-1)$-ortho-symmetric} or \emph{$n$-precluster tilting} if $\add(_\Lambda M)\subseteq {}^{\perp_{n-1}}M=M^{\perp_{n-1}}$. The algebra $A$ is $n$-Auslander if and only if there is an algebra $\Lambda$ and a maximal $(n-1)$-orthogonal $\Lambda$-module ${_\Lambda}N$ such that $A=\End_\Lambda(N)$ by \cite[Proposition 2.4.1]{Iyama07}, and is $n$-minimal Auslander-Gorenstein if and only if  there is an algebra $\Lambda$ and an $(n-1)$-ortho-symmetric generator-cogenerator ${_\Lambda}N$ such that $A=\End_\Lambda(N)$ by \cite[Theorem 4.5]{IS} or \cite[Corollary 3.18]{CK}. Moreover, by \cite[Proposition 4.1]{IS}, if $A$ is $n$-minimal Auslander-Gorenstein, then either $A$ is self-injective or $\id({_A}A)=n+1=\dm(A)$. In the latter case, we  have $\id(A{_A})=n+1=\dm(A)$ and thus  $A$ is $(n+1)$-Gorenstein.

An $A$-module $M$ is said to be \emph{$m$-rigid} if $\Ext_A^i(M,M)=0$ for all $1\leq i\leq m$. Over symmetric algebras, ortho-symmetric modules have been characterized as  follows.

\begin{Lem}{\rm \cite[Corollary 5.4]{CK}}\label{OSS}
Let $A$ be a symmetric algebra and $N$ a basic $A$-module without any nonzero projective direct summands. For a natural number $m$. the $A$-module $A\oplus N$ is $m$-ortho-symmetric if and only if $N$ is $m$-rigid and $\Omega_A^{m+2}(N)\cong N$.
\end{Lem}

\section{Mirror-reflective algebras}

In this section, we introduce (reduced) mirror-reflective  algebras and prove that these algebras can be linked by recollements of their derived module categories (see Proposition \ref{Idem-rec}).

\subsection{Mirror-reflective algebras and their derived recollements\label{sect3.1}}
Throughout this section, assume that $A$ is an algebra over a commutative ring $k$. Let $M$ be an $A$-$A$-bimodule and $\alpha: {}_AM\otimes_AM\ra M$ be a homomorphism of $A$-$A$-bimodules, such that the associative law holds
$$ (\heartsuit)\qquad \big((x\otimes y)\alpha \otimes z\big)\alpha= \big(x\otimes (y \otimes z)\alpha\big)\alpha \; \mbox{ for } x,y,z\in M.$$
We define a multiplication on the underlying abelian group $A\oplus M$ by setting
$$ (a, m)\cdot (b,n):= (ab, an+mb+(m\otimes n)\alpha)\; \mbox{ for } a,b\in A, \, m,n\in M. $$
Then $A\oplus M$ becomes an associative algebra with the identity $(1,0)$, denoted by $R(A,M,\alpha)$. In the following, we identify $A$ with $(A,0)$, and $M$ with $(0,M)$ in $R(A,M,\alpha)$. Thus $A$ is a subalgebra of $R(A,M,\alpha)$ with the same identity, and $M$ is an ideal of $R(A,M,\alpha)$ such that $R(A,M,\alpha)/M\simeq A$.

\smallskip
Now, we consider a special case of the above construction. Let $e=e^2\in A$, $\Lambda:=eAe\;\;\mbox{and}\;\; Z(\Lambda)$ be the center of $\Lambda$. For $\lambda\in Z(\Lambda)$,  let $\omega_\lambda$ be the composite of the natural maps:
$$(Ae\otimes_\Lambda eA)\otimes_A(Ae\otimes_\Lambda eA)\lraf{\simeq}Ae\otimes_\Lambda (eA\otimes_AAe)\otimes_{\Lambda} eA \lraf{\simeq} Ae\otimes_\Lambda \Lambda\otimes_\Lambda eA\lraf{{\rm Id}\otimes(\cdot\lambda)\otimes {\rm Id}}
Ae\otimes_\Lambda \Lambda\otimes_\Lambda eA\ra Ae\otimes_\Lambda eA,
$$where $(\cdot\lambda):\Lambda\to \Lambda$ is the multiplication map by $\lambda$.
Then $\omega_{\lambda}$ satisfies the associative law $(\heartsuit)$.

Let $R(A,e,\lambda):=R(A, Ae\otimes_{}eA,\omega_{\lambda})$. Then the elements of $R(A,e, \lambda)$
are of the form
$$a+\sum_{i=1}^n a_ie\otimes eb_i \;\;\mbox{for}\;\;a, a_i, b_i\in A, 1\leq i\leq n\in\mathbb{N}.$$
The multiplication, denoted by $\ast$, is explicitly given by
$$ (a + be\otimes ec)\ast (a' + b'e\otimes ec'):= aa' + (ab'e\otimes ec'+ be\otimes eca'+ becb'e\otimes\lambda ec')$$
for $a,b,c,a',b',c'\in A$, and can be extended linearly to elements of general form.
Particularly, $$(\diamondsuit)\;\quad (ae\otimes eb)\ast(a'e\otimes eb')=aeba'e\lambda\otimes eb'=ae\otimes \lambda eba'e b'.$$

\begin{Def}\label{RMSA}
The algebra $R(A,e,\lambda)$ defined above is called the \emph{mirror-reflective algebra} of $A$ at level $(e,\lambda)$. The algebra $R(A,e,e)$ is then called the \emph{mirror-reflective algebra} of $A$ at $e$, denoted by $R(A,e)$.

The algebra $S(A,e,e):=(1-e)A(1-e)\oplus Ae\otimes_\Lambda eA$ with the multiplication induced from the one of $R(A,e,e)$ is called the \emph{reduced mirror-reflective  algebra} of $A$ at $e$, denoted by $S(A,e)$.
\end{Def}

Compared with $R(A,e)$, $S(A,e)$ has a fewer number of simple modules. So it is termed the reduced mirror-reflective  algebra. The following lemma is obvious.

\begin{Lem}\label{rc}
$(1)$ There is an algebra isomorphism $R(A,e,\lambda)/(Ae\otimes_\Lambda eA)\simeq A$.

$(2)$ If $\mu\in Z(\Lambda)$ is an invertible element, then $R(A,e,\lambda)\simeq R(A,e,\lambda\mu)$ as algebras.
\end{Lem}

\medskip
For simplicity, let $R:=R(A,e),\;\;S:=S(A,e)\;\;\mbox{and}\;\;\bar{e}:=e\otimes e\in R.$
Then $\bar{e}=\bar{e}^2$, $e\bar{e}=\bar{e}=\bar{e}e$, and $\{\bar{e}, e-\bar{e}, 1-e\}$ is a set of orthogonal idempotents in $R$. Now, we define
$$\pi_1: R\lra A,\;\;a+\sum_{i=1}^n a_i\bar{e}b_i\mapsto a, \mbox{ and } \;
\pi_2: R\lra A,\;\;a+\sum_{i=1}^n a_i\bar{e}b_i\mapsto a+\sum_{i=1}^n a_ieb_i$$
for $a, a_i, b_i\in A$ and $1\leq i\leq n.$
Then $\pi_1$ and $\pi_2$ are surjective homomorphisms of algebras.
Let
$$I:=\Ker(\pi_1), \;\;J:=\Ker(\pi_2)\;\;\mbox{and}\;\;e_0:=(1-e)+\overline{e}.$$

\begin{Lem}\label{Property}
$(1)$ $I=R\bar{e}R$, $J=R(e-\bar{e})R$, $IJ=0=JI$, $I+J=ReR$ and $S=e_0Re_0$.

$(2)$ As an $A$-$A$-bimodule, ${}_AR_A$ has two decompositions:
$R=A\oplus I=A\oplus J.$

$(3)$ There is an algebra isomorphism
$\phi: R\ra R,\;\; a+\sum_{i=1}^n a_i\bar{e}b_i\mapsto a+\sum_{i=1}^n a_i(e-\bar{e})b_i.$
Moreover, $\phi^2={\rm Id}_R$ and the restriction of $\phi$ to $I$ induces an isomorphism $I\to J$ of $A$-$A$-bimodules.

$(4)$ $\pi_2=\phi\pi_1$ and both $\pi_1$ and $\pi_2$ induce surjective homomorphisms of algebras
$$\pi_1': S\lra (1-e)A(1-e)\quad\mbox{and}\quad\pi_2':S\lra A,$$
respectively.  Moreover, $\Ker(\pi_1')=I$ and $\Ker(\pi_2')=(1-e)J(1-e)=J\cap S$.
\end{Lem}

{\it Proof.} $(1)$ Clearly, $I=Ae\otimes_\Lambda eA=A\bar{e}A= R\bar{e}R$.  Since $(e-\bar{e})\pi_2=0$, we have $e-\bar{e}\in\Ker(\pi_2)=J$ and $R(e-\bar{e})R\subseteq J$.
Conversely, if $r:=a+\sum_{i=1}^n a_i\bar{e}b_i\in J$, then $a+\sum_{i=1}^n a_ieb_i=(r)\pi_2=0$, that is, $a=-\sum_{i=1}^n a_ieb_i$.
Consequently, $r=-\sum_{i=1}^n a_ieb_i+\sum_{i=1}^n a_i\bar{e}b_i=-\sum_{i=1}^n a_i(e-\bar{e})b_i\in R(e-\bar{e})R$. Thus $J=R(e-\bar{e})R=A(e-\bar{e})A$.
Note that $I+J=R\bar{e}R+R(e-\bar{e})R=ReR$. For any  $x,y, x', y'\in A$, since $(x\bar{e}y)(x'(e-\bar{e})y')=x\bar{e}yx'ey'-x\bar{e}yx'ey'=0$, we have $IJ=0$.
Similarly, $(x'(e-\bar{e})y')(x\bar{e}y)=0$, and therefore $JI=0$. Since $I$ is an ideal of $R$ and $IJ=JI=0$, it follows that $S=e_0Re_0$.

$(2)$ $R$ contains $A$ as a subalgebra with the same identity, and the composite of the inclusion $A\subseteq R$ with $\pi_i$ for $i=1,2$, is the identity map of $A$. Thus $(2)$ follows.

$(3)$ By $(2)$, $I\simeq R/A\simeq J$ as $A$-$A$-bimodules. In fact, the isomorphism from
$I$ to $J$ is given by
$$
\varphi': \; I\lra J,\;\;\sum_{i=1}^n a_i\bar{e}b_i\mapsto \sum_{i=1}^n a_i (e-\bar{e})b_i.
$$
Then the map
$\phi: R=A\oplus I\to R=A\oplus J$ induced from $\varphi'$ is a well-defined  isomorphism of $A$-$A$-bimodules. Moreover, $\phi$ preserves the multiplication of $R$ and $\phi^2={\rm Id}_R$. Thus $\phi$ is an automorphism of algebras.

$(4)$ The first equality in $(4)$ follows from the definitions of $\phi$, $\pi_1$ and $\pi_2$. To see other statements in $(4)$, we apply the left and right multiplications by $e_0$ to $\pi_1$ and $\pi_2$, and then use $(1)$.   $\square$

\medskip
Under some conditions, the associated ideals $I$ and $J$
are related by annihilators of modules. Recall that the \emph{annihilator} of an $R$-module $M$
is defined as ${\rm Ann}_{R}(M):=\{r\in R\mid rM=0\}$, which is an ideal of $R$.

\begin{Lem}\label{ring homo}
$(1)$ If $eA_{A}$ is a faithful right module, then $J={\rm Ann}_{R\opp}(I)$. Dually, if $_{A}Ae$ is faithful, then $J={\rm Ann}_R(I)$.

$(2)$ The map $\pi_2$ induces isomorphisms of abelian groups:
$$R\bar{e}\lraf{\simeq} Ae,\;\; \bar{e}R\lraf{\simeq}eA\;\;\mbox{and}\;\; \bar{e}R\bar{e}\lraf{\simeq}eAe.$$
The map $\pi_2':S\to A$ in Lemma \ref{Property}(4) induces isomorphisms of abelian groups:
$$S\bar{e}\lraf{\simeq} Ae,\;\; \bar{e}S\lraf{\simeq}eA\;\;\mbox{and}\;\; \bar{e}S\bar{e}\lraf{\simeq}eAe.$$

$(3)$ The map $\pi_1$ induces isomorphisms of abelian groups:
$$R(e-\bar{e})\lraf{\simeq} Ae,\;\; (e-\bar{e})R\lraf{\simeq}eA\;\;\mbox{and}\;\; (e-\bar{e})R(e-\bar{e})\lraf{\simeq}eAe.$$
\end{Lem}

{\it Proof.} $(1)$ Since $IJ=0$ by Lemma \ref{Property}(1), $J\subseteq{\rm Ann}_{R\opp}(I)$. Let $x:=a+\sum_{i=1}^n a_i\bar{e}b_i\in {\rm Ann}_{R\opp}(I)$ with $a, a_i,b_i\in A$ and $1\leq i\leq n\in\mathbb{N}$. Since $J=\Ker(\pi_2)$, it suffices to show $y:=(x)\pi_2=0$.  In fact, by $Ix=0$, we have $0=(Ix)\pi_2=(I)\pi_2 y=AeAy$. This implies $eAy=0$. In other words, $y\in{\rm Ann}_{A\opp}(eA)$. Since $eA_{A}$ is faithful, ${\rm Ann}_{A\opp}(eA)=0$. Thus $y=0$. This shows $J={\rm Ann}_{R\opp}(I)$. Similarly, we show the second identity.

$(2)$ Since $(\bar{e})\pi_2=e$, the restriction $f_2: R\bar{e}\to Ae$ of $\pi_2$ to $R\bar{e}$ is surjective.
As $\Ker(f_2)=R\bar{e}\cap J\subseteq JI=0$ by Lemma \ref{Property}(1),  $f_2$ is an isomorphism.
Dually, the restriction $\bar{e}R\to eA$ of $\pi_2$ to $\bar{e}R$ is also an isomorphism.
Consequently, $\pi_2$ induces an algebra isomorphism from $\bar{e}R\bar{e}$ to $eAe$.

Since $IJ=JI=0$ by Lemma \ref{Property}(1),  we have $S\overline{e}=R\overline{e}$ and $\overline{e}S=\overline{e}R$. Clearly,  $\overline{e}S\overline{e}=\overline{e}R\overline{e}$. Thus the second statement in $(2)$ holds.

$(3)$ This follows from $(2)$ and Lemma \ref{Property}(3)-(4). $\square$

Consequently, Lemma \ref{Property}(1) and Lemma \ref{ring homo}(2) imply $\#(R)=\#(A)+\#(eAe)$.

\begin{Prop}\label{Faithful}
Let $A_2:=\End_R(R\oplus R/I)$ and $B_2:=\End_{S}(S\oplus S/I)$. Suppose that the right $A$-module $eA_A$ is faithful. Then the following hold true.

$(1)$  There are standard recollements of derived module categories

$$\xymatrix{\mathscr{D}(A)\ar[r]&\mathscr{D}(A_2)\ar[r]\ar@/^1.2pc/[l]\ar@/_1.2pc/[l]
&\mathscr{D}(A),\ar@/^1.2pc/[l]\ar@/_1.2pc/[l]}\vspace{0.2cm}\quad \xymatrix{\mathscr{D}(A)\ar[r]&\mathscr{D}(B_2)\ar[r]\ar@/^1.2pc/[l]\ar@/_1.2pc/[l]
&\mathscr{D}((1-e)A(1-e))\ar@/^1.2pc/[l]\ar@/_1.2pc/[l]}\vspace{0.2cm}$$
induced by finitely generated, right-projective idempotent ideals of $A_2$ and $B_2$, respectively.

$(2)$ $\sd(A_2)\geq 2\sd(A)+1$ and $\sd(B_2)\geq \sd(A)+\sd((1-e)A(1-e))+1$.

$(3)$ $\gd(A\opp)\leq \gd(A_2\opp)\leq 2\gd(A\opp)+2,$ $\fd(A\opp)\leq \fd(A_2\opp)\leq 2\fd(A\opp)+2.$
\end{Prop}

{\it Proof.}  $(1)$ Since $eA_A$ is faithful, $J={\rm Ann}_{R\opp}(I)$ by Lemma \ref{ring homo}(1). Note that $I$ is an ideal of $S$ and
${\rm Ann}_{S\opp}(I)=S\cap{\rm Ann}_{R\opp}(I)=S\cap J$. By Lemma \ref{Property}(4), there are algebra isomorphisms
$A\simeq R/I\simeq R/J\simeq S/(S\cap J)$ and $S/I\simeq (1-e)A(1-e)$. Now, Proposition \ref{Faithful} follows from Lemma \ref{Projective}(1).

$(2)$ This follows from $(1)$ and Proposition \ref{upper}(3).

$(3)$ This is a consequence of Proposition \ref{Faithful}(1) and
\cite[Corollary 3.12 and Theorem 3.17]{xc7}. $\square$

\medskip
Now, we consider $n$-idempotent and strong idempotent ideals of mirror-reflective  algebras.

\begin{Prop}\label{Idem-rec}
$(1)$ The ideals $I$ and $J$ of $R$ are $2$-idempotent.

$(2)$  Let $n\ge 1$ be an integer. Then $I$ is $(n+2)$-idempotent if and only if so is $J$ if and only if $\Tor_i^{eAe}(Ae, eA)=0$ for all $1\leq i\leq n$.

$(3)$ If $\Tor_i^{eAe}(Ae, eA)=0$ for all $i\geq 1$, then there are standard recollements of derived module categories induced by $I:=R\bar{e}R$:
\smallskip
$$\xymatrix{\mathscr{D}(A)\ar[r]&\mathscr{D}(R)\ar[r]\ar@/^1.2pc/[l]\ar@/_1.2pc/[l]
&\mathscr{D}(eAe) \; \; \; \mbox{ and}\ar@/^1.2pc/[l]\ar@/_1.2pc/[l]}\vspace{0.2cm}\quad \xymatrix{\mathscr{D}((1-e)A(1-e))\ar[r]&\mathscr{D}(S)\ar[r]\ar@/^1.2pc/[l]\ar@/_1.2pc/[l]
&\mathscr{D}(eAe)\ar@/^1.2pc/[l]\ar@/_1.2pc/[l]}\vspace{0.2cm}.$$
\end{Prop}

{\it Proof.} $(1)$ There is a commutative diagram
$$
\xymatrix{R\overline{e}\otimes_{\overline{e}R\overline{e}}\overline{e}R\ar[r]^-{\mu}\ar[d]_-{\pi_2\otimes\pi_2}& R\overline{e}R=Ae\otimes_{eAe}eA\ar[d]_-{\pi_2}\\
Ae\otimes_{eAe}eA\ar[r]^-{\mu'}& AeA}
$$
where $\mu$ and $\mu'$ are given by multiplications. By Lemma \ref{ring homo}(2), $\pi_2\otimes\pi_2$ is an isomorphism.
Note that the composition of the inverse of $\pi_2\otimes\pi_2$ with $\mu$ is the identity of $Ae\otimes_{eAe}eA$. Thus $\mu$
is an isomorphism. This shows that $I$ is $2$-idempotent by Lemma \ref{strong idempotent}(1). Similarly, we can show that $J$ is $2$-idempotent by using the idempotent $e-\overline{e}$ and the algebra homomorphism $\pi_1$.

$(2)$ By Lemma \ref{Property}(3), $I$ is $(n+2)$-idempotent if and only if so is $J$. Since $I$ is $2$-idempotent by $(1)$, it follows from
Lemma \ref{strong idempotent}(1) that $I$ is $(n+2)$-idempotent if and only if $\Tor_i^{\overline{e}R\overline{e}}(R\overline{e}, \overline{e}R)=0$
for $1\leq i\leq n$. By Lemma \ref{ring homo}(2),  $\pi_2$ induces isomorphisms of abelian groups
$\Tor_i^{\overline{e}R\overline{e}}(R\overline{e}, \overline{e}R)\simeq \Tor_i^{eAe}(Ae, eA)$ for all $i\in\mathbb{N}$. Thus $I$ is $(n+2)$-idempotent if and only if
$\Tor_i^{eAe}(Ae, eA)=0$ for $1\leq i\leq n$.

$(3)$ By $(2)$, $I$ is a strong idempotent ideal of $R$ if and only if $\Tor_i^{eAe}(Ae, eA)=0$ for all $i\geq 1$. According to Corollary \ref{HID}(1), if $I$ is a strong idempotent ideal of $R$, then $e_0Ie_0$ is a strong idempotent ideal of $S$. By Lemma \ref{Property} and Lemma \ref{ring homo}(2), $e_0Ie_0=I$, $S/I\simeq (1-e)A(1-e)$, $R/I\simeq A$ and $\overline{e}R\overline{e}\simeq eAe\simeq \overline{e}S\overline{e}$.
Thus the recollements in $(3)$ exist. $\square$

\medskip
To discuss the decomposition of $R$ as an algebra and to lift algebra homomorphisms, we show the following result.
For a homomorphism $\alpha: A\to \Gamma$ of algebras, denote by $\Hom_{\alpha\mbox{-}\rm Alg}(R,\Gamma)$ the set of all algebra homomorphisms
$\beta: R\to \Gamma$ such that the restriction of $\beta$ to $A$ coincides with $\alpha$.

\begin{Lem}\label{generator}
$(1)$ If $u=u^2\in A$ such that $\add({_A}Au)=\add({_A}Ae)$, then $R\simeq R(A,u,u)$ as algebras.

$(2)$ If ${_A}Ae$ is a generator, then $R\simeq A\times A$ as algebras.

$(3)$ If $\Gamma$ is an algebra and $\alpha:A\to \Gamma$ is an algebra homomorphism, then there is a bijection
$$
\Hom_{\alpha\mbox{-}\rm Alg}(R,\Gamma)\lraf{\simeq}\{x\in (e)\alpha \; \Gamma \;(e)\alpha\mid x^2=x,\; (c)\alpha\, x=x \, (c)\alpha\;\;\mbox{ for }\;\; c\in\Lambda\},\;\;
\overline{\alpha}\mapsto (\overline{e})\overline{\alpha}.
$$
\end{Lem}

{\it Proof.} $(1)$ Let $U:=uAu$. We keep the notation in the proof of Lemma \ref{strong idempotent} and identify $\Hom_A(Au,-): A\Modcat\to U\Modcat$ with the functor $u\cdot: A\Modcat\to U\Modcat$, given by the left multiplication of $u$. Let $\mu: Au\otimes_U u(-)\to {\rm Id}$ be the counit adjunction of the adjoint pair $(Au\otimes_U-, u\cdot)$. Then, for an $A$-module $X$, the map $\mu_X$ is an isomorphism if and only if $X\in \mathbf{P}_1(Au)$. Applying $Ae\otimes_{\Lambda}-$ to a projective presentation of ${_\Lambda}eA$, we obtain an exact sequence $P_1\to P_0\to Ae\otimes_{\Lambda}eA\to 0$ of $A$-modules with $P_1, P_0\in\Add(Ae)$. This shows $Ae\otimes_\Lambda eA\in\mathbf{P}_1(Ae)$. Due to $\add({_A}Au)=\add({_A}Ae)$, we have $Ae\otimes_\Lambda eA\in\mathbf{P}_1(Au)$, and therefore
$\mu_{Ae\otimes_\Lambda eA}: Au\otimes_U u(Ae\otimes_\Lambda eA)\to Ae\otimes_\Lambda eA$
is an isomorphism of $A$-$A$-modules. Since the multiplication map $\rho: Ae\otimes_\Lambda eA\to A, ae\otimes eb\mapsto aeb$ for $a, b\in A$, satisfies $e\Ker(\rho)=0=e\Coker(\rho)$, it follows from $\add({_A}Au)=\add({_A}Ae)$ that $u\Ker(\rho)=0=u\Coker(\rho)$. Then $u\rho: u(Ae\otimes_\Lambda eA)\to uA$ is an isomorphism of $U$-$A$-bimodules, and $u\rho u: u(Ae\otimes_\Lambda eA)u\to uAu$ is an isomorphism of $U$-$U$-bimodules. Consequently, there is an isomorphism of  $A$-$A$-bimodules $$ Au\otimes_U u\rho: \;\;Au\otimes_U u(Ae\otimes_\Lambda eA)\lraf{\simeq} Au\otimes_U uA.$$
Thus  $\psi:=(Au\otimes_U u\rho)^{-1}\mu_{Ae\otimes_\Lambda eA}: Au\otimes_U uA\to Ae\otimes_\Lambda eA$ is an isomorphism of $A$-$A$-bimodules. In fact, if $x_i\in uAe$ and $y_i\in eAu$ with $1\leq i\leq n$ such that $\sum_{i=1}^nx_iy_i=u$, then $(a(u\otimes u)b)\psi=a(\sum_{i=1}^nx_i\otimes y_i)b.$
This induces an isomorphism of $A$-$A$-bimodules: $$({\rm Id}_A, \psi): R(A,u,u)=A\oplus Au\otimes_U uA\lra R=A\oplus Ae\otimes_\Lambda eA,$$ $$ (a, x\otimes y)\mapsto (a, (x\otimes y)\psi) \mbox{ for } a\in A, x\in Au, y\in uA. $$
A verification shows that this is an algebra isomorphism.

$(2)$ Suppose that ${_A}Ae$ is a generator. Then $\add(_AAe)=\add(_AA)$. Let $B:=R(A,1,1)$. By $(1)$, $R\simeq B$ as algebras. Now, identifying $A\otimes_AA$ with $A$, we then get $B=A\oplus A$ with the multiplication given by $$(a_1,a_2)(b_1,b_2):=(a_1b_1, a_1b_2+a_2b_1+a_2b_2) \; \mbox{ for } a_1,a_2, b_1, b_2\in A.$$
Clearly, $(1,0)$ is the identity of $B$ and $(1,-1)$ is a central idempotent of $B$. Thus
the map $B\to A\times A$, $(a_1,a_2)\mapsto (a_1, a_1+ a_2)$, is an algebra isomorphism. Thus $R\simeq B\simeq A\times A$ as algebras.

$(3)$  Note that $\Gamma$ can be regarded as an $A$-$A$-bimodule via $\alpha$ and that any $A$-$A$-bimodule can be considered as a module over the enveloping algebra $A^e:=A\otimes_kA^{\opp}$. Define $F=Ae\otimes_\Lambda-\otimes_\Lambda eA: \Lambda^{e}\Modcat\to A^{e}\Modcat$ and $G=e(-)e:  A^{e}\Modcat\to \Lambda^{e}\Modcat.$ Then there are isomorphisms of $k$-modules
$$\Hom_{A^{\rm e}}(Ae\otimes_\Lambda eA, \Gamma)\simeq\Hom_{A^e}(F(\Lambda), \Gamma)\simeq\Hom_{\Lambda^e}(\Lambda, G(\Gamma))=\Hom_{\Lambda^e}\big(\Lambda, (e)\alpha\, \Gamma\, (e)\alpha\big)$$
$$\quad\quad\quad\quad \quad=\{y\in (e)\alpha\; \Gamma \,(e)\alpha \mid (c)\alpha\,y =y\, (c)\alpha\;\;\mbox{for any}\;\; c\in\Lambda\}=:\Gamma'.$$
Let $\overline{\alpha}\in \Hom_{\alpha\mbox{-}\rm Alg}(R,\Gamma)$ and $x:=(\overline{e})\overline{\alpha}\in \Gamma$.
Since the restriction of $\overline{\alpha}$ to $A$ equals $\alpha$, the restriction of $\overline{\alpha}$ to $Ae\otimes_\Lambda eA$ is an homomorphism of $A$-$A$-bimodules. By $\overline{e}^2=\overline{e}$, we have $x^2=x\in \Gamma'$ and $(ae\otimes eb)\overline{\alpha}=(a)\alpha \, x\, (b)\alpha$ for any $a,b\in A$. This means that $\overline{\alpha}$ is determined by $\alpha$ and $x$.

Conversely, let $y\in \Gamma'$ and let $h: Ae\otimes_\Lambda eA\to \Gamma$ be the homomorphism of $A$-$A$-bimodules sending $ae\otimes eb$ to $(a)\alpha \; y\; (b)\alpha$.
Define $\overline{h}:=(\alpha, h): R\to \Gamma$. Then $\overline{h}$ is an algebra homomorphism if and only if
$((ae\otimes eb)\ast(a'e\otimes eb'))h=(ae\otimes eb)h(a'e\otimes eb')h$ for any $a, a', b, b'\in A$ if and only if $y(ba')\alpha \; y=(eba')\alpha \; y$ for any $b,a'\in A$.
Now, suppose $y^2=y$. Since $\alpha$ is an algebra homomorphism and $(e)\alpha \; y=y=y \; (e)\alpha$,
we see that $(eba')\alpha \; y=(eba'e)\alpha\; y=(eba'e)\alpha \; y^2=y\; (eba'e)\alpha y=y (ba')\alpha \, y.$ Thus $\overline{h}$ is an algebra homomorphism.
As $y=(\overline{e})h$, the bijection in $(3)$ is clear.
$\square$

\begin{Prop} \label{Indecomposable}
Let $A$ be an indecomposable algebra. Then

$(1)$ $R$ is a decomposable algebra if and only if ${_A}Ae$ is a generator.
In this case, $R\simeq A\times A$ as algebras.

$(2)$ If $\add(Ae)\cap\add(A(1-e))=0$ and $(1-e)A(1-e)$ is an indecomposable algebra, then $S$ is an indecomposable algebra.
\end{Prop}

{\it Proof.}
$(1)$ If ${_A}Ae$ is a generator, then $R\simeq R(A,1,1)\simeq A\times A$ as algebras by Lemma \ref{generator}(2), and therefore $R$ is decomposable.
Conversely, assume that $R$ is a decomposable algebra. Then there is an idempotent $z$ in the center $Z(R)$ of $R$ such that $z\neq 0,1$.
Since $\pi_1:R\to A$ is a surjective homomorphism of algebras, it restricts to an algebra homomorphism $Z(R)\to Z(A)$. This implies $(z)\pi_1\in Z(A)$.
Since $A$ is indecomposable, $(z)\pi_1=0$ or $1$. If  $(z)\pi_1=0$, then $z\in I=\Ker(\pi_1)$. If $(z)\pi_1=1$, then $1-z\in I$.
So, without loss of generality, we can assume $z\in I$. Similarly, $z\in J$ or $1-z\in J$ by $\pi_2$.
If $z\in J$, then $z=z^2\in IJ=0$ by Lemma \ref{Property}(1), a contradiction. Thus $1-z\in J$, and  $1=z+(1-z)\in I+J=ReR$ by Lemma \ref{Property}(1). This shows $ReR=R$ and implies $AeA=A$ by $\pi_1$. Hence ${_A}Ae$ is a generator.

$(2)$ Let $J_1:=S\cap J$.  In the proof of $(1)$, we replace $\pi_1$ and $\pi_2$ with $\pi_1': S\to (1-e)A(1-e)$ and $\pi_2': S\to A$ (see Lemma \ref{Property}(4)), respectively, and show similarly that if $(1-e)A(1-e)$ is indecomposable and $S$ is decomposable, then $S=I+J_1$. In this case, the equality $A=AeA$ still holds because $\pi_2'$ is surjective with $\Ker(\pi_2')=J_1$ and $(\overline{e})\pi_2'=e$. Consequently, ${_A}Ae$ is a generator, and therefore the assumption $\add(Ae)\cap\add(A(1-e))=0$ forces $e=1$. Thus $S=I\simeq A$ as algebras. This contradicts to $A$ being indecomposable. $\square$

\subsection{Examples of mirror-reflective algebras: quivers with relations} \label{MQA}

In this subsection, we describe explicitly the mirror-reflective algebras for algebras presented by quivers with relations. This explains the terminology ``mirror-reflective algebras" (see Example \ref{Exmaple1} below).

Let $Q:=(Q_0,Q_1)$ be a quiver with the vertex set $Q_0$ and arrow set $Q_1$. For an arrow $\alpha:i\to j$, we denote by $s(\alpha)$ and $t(\alpha)$ the starting vertex $i$ and the terminal vertex $j$, respectively. Composition of an arrow $\alpha: i\to j$ with an arrow $\beta: j\to m$ is written as $\alpha\beta$.
A \emph{path} of length $n\ge 0$ in $Q$ is a sequence $p:=\alpha_1\cdots\alpha_n$ of $n$ arrows $\alpha_i$ in $Q_1$ such that $t(\alpha_i)=s(\alpha_{i+1})$ for $1\leq i<n\in\mathbb{N}$. Set $s(p)=s(\alpha_1)$ and $t(p)=t(\alpha_n)$. In case of $n=0$, we understand the trivial path as an vertex $i\in Q_0$, denote by $e_i$ and set $s(e_i)=i=t(e_i)$. We write $\mathscr{P}(Q)$ for the set of all paths of finite length in $Q$.
For a field $k$, we write $kQ$ for the path algebra of $Q$ over $k$. Clearly, it has $\mathscr{P}(Q)$ as a $k$-basis.

A \emph{relation} $\sigma$ on $Q$ over $k$ is a $k$-linear combination of paths $p_i$ of length at least $2$. We may assume that all paths in a relation have the same starting vertex and terminal vertex, and define $s(\sigma)=s(p_i)$ and $t(\sigma)=t(p_i)$. If $\rho=\{\sigma_i\}_{i\in T}$ is a set of relations on $Q$ over $k$ with $T$ an index set, the pair $(Q,\rho)$ is called a \emph{quiver with relations} over $k$. In this case, we have a $k$-algebra $k(Q,\rho):=kQ/\langle\rho\rangle$, the quotient algebra of the path algebra $kQ$ modulo the ideal  $\langle\rho\rangle$ generated by the relations $\sigma_i, i\in T$.


\begin{Lem}\label{quiver}
Let $B$ be a $k$-algebra, $\{f_i\mid i\in Q_0\}$ a set of orthogonal idempotents in $B$ with $1_B=\sum_{i\in Q_0}f_i$, and $\{f_\alpha\mid \alpha\in Q_1\}$ a set of elements in $B$. If $f_{s(\alpha)}f_\alpha=f_\alpha=f_\alpha f_{t(\alpha)}$ for  $\alpha\in Q_1$, then there is a unique algebra homomorphism $f: kQ\to B$ which sends $e_i\mapsto f_i$ and $\alpha\mapsto f_\alpha$.
\end{Lem}

Let $Q':=(Q'_{0}, Q'_{1})$ be a full subquiver of $Q$, that is, $Q'_{0}\subseteq Q_0$ and $Q'_{1}=\{\alpha\in Q_1\mid s(\alpha), t(\alpha)\in Q'_{0}\}$.
Define $$A:=k(Q,\rho), \;\;V_0:=Q_0\setminus Q'_{0} \; \; \mbox{and}\;\; e:=\sum_{i\in V_0}e_i\in A.$$
We shall describe the quiver and relations for the mirror-reflective algebra $R(A,e)$ explicitly.

Let $\overline{Q}$ be a copy of the quiver $Q$, say $\overline{Q}_0=\{\bar{i}\mid i\in Q_0\}$ and $\overline{Q}_1=\{\bar{\alpha}\mid \alpha\in Q_1\}$, with $s(\bar{\alpha})=\bar{i}$ and $t(\bar{\alpha})=\bar{j}$ if $s(\alpha)=i$ and $t(\alpha)=j$. Consider $Q'$ as a full subquiver of $\overline{Q}$ by identifying $\bar{i}$ with $i$ for $i\in Q'_0$, and $\bar{\alpha}$ with $\alpha$ for $\alpha\in Q'_1$. So $Q_0\cap \overline{Q}_0=Q'_0$ and $Q_1\cap \overline{Q}_1 = Q'_1$. Let $\Delta:=(\Delta_0,\Delta_1)$ be the pullback of $Q$ and $\overline{Q}$ over $Q'$, that is, $$\Delta_0:= Q_0 \dot{\cup} (\overline{Q_{0}}\setminus Q'_0)\quad\mbox{and}\quad \Delta_1:=Q_1 \dot{\cup} (\overline{Q_1}\setminus Q'_{1}).$$
We define a map $(-)^+:\{e_i\mid i\in Q_0\}\cup Q_1\to k\Delta$ by
$$e_i^+:=\left\{
\begin{array}{ll}
e_i, & \hbox{$i\in Q'_{0}$},\\
e_i+e_{\overline{i}}, & \hbox{$i\in V_{0}$},
\end{array}
\right. \quad
\alpha^+:=\left\{
\begin{array}{ll}
\alpha, & \hbox{$\alpha\in Q'_{1}$,} \\
\alpha+\overline{\alpha}, & \hbox{$\alpha\in Q_1\setminus Q'_{1}$.}
\end{array}
\right.
$$
Since $e_{s(\alpha)}^+\alpha^+=\alpha^+=\alpha^+e_{t(\alpha)}^+$ for any $\alpha\in Q_1$, it follows from Lemma \ref{quiver}
that  $(-)^+$ can be extended to an algebra homomorphism $$(-)^+:\quad kQ\lra k\Delta,\quad p\mapsto p^+:=\alpha_1^+\cdots\alpha_n^+ \mbox{ for } p=\alpha_1\cdots\alpha_n\in \mathscr{P}(Q).$$

Given a relation $ \sigma:=\sum_{i=1}^na_ip_i$ on $Q$ with $a_i\in k$, $p_i\in\mathscr{P}(Q)$ for $1\leq i\leq n\in\mathbb{N}$,
and $s(\sigma), t(\sigma)\in Q'_{0}$, we define
$$
\sigma_+:=\sum_{1\le j\le n,\; p_j\in \mathscr{P}(Q')}a_jp_j \; +\sum_{1\le i\le n,\; p_i\notin\mathscr{P}(Q')}a_i(p_i+\overline{p_i})=\sigma \; + \sum_{1\le i\le n,\; p_i\notin\mathscr{P}(Q')}a_i\overline{p_i}.
$$
Now, let $\psi :=\psi_1\cup \psi_2\cup\psi_3\cup\psi_4$ with
$$\begin{array}{ll} \psi_1 :=&\{\overline{a}p b, a p\overline{b} \mid a, b\in Q_1, s(a), t(b)\in V_{0}, p\in\mathscr{P}(Q'), apb\in\mathscr{P}(Q)\},\\
\psi_2 :=&\{\sigma\in\rho\mid s(\sigma)\in V_{0}\;\;\mbox{or}\;\;t(\sigma)\in V_{0}\},\\
\psi_3 :=&\{\overline{\sigma}\mid\sigma\in\psi_2\},\; \mbox{and}\\
\psi_4:=&\{\sigma_+\mid \sigma\in\rho, s(\sigma), t(\sigma)\in Q'_{0}\}.\\
\end{array}$$
Then $\psi$ is a set of relations on $\Delta$ over $k$, and we consider the $k$-algebra $k(\Delta, \psi)$.

\begin{Prop}\label{Quiver algebra}
$(1)$ The homomorphism $(-)^+:kQ\to k\Delta$ of algebras is injective and induces an injective homomorphism $\mu: A\to k(\Delta, \psi)$ of algebras.

$(2)$ There exists an isomorphism $\theta: R(A,e)\lraf{\simeq} k(\Delta, \psi)$ of algebras such that $(e_i\otimes e_i)\theta=e_{\overline{i}}$ for $i\in V_0$, and the restriction of $\theta$ to $A$
coincides with $\mu$ in \emph{(1)}.
\end{Prop}

{\it Proof.} $(1)$ For $\mathcal{U}\subseteq k\Delta$, we denote by $\langle\mathcal{U}\rangle$
the ideal of $k\Delta$ generated by $\mathcal{U}$. Let $E:=\{e_{\overline{i}}\mid i\in V_{0}\}$. Then $k\Delta/\langle E\rangle\lraf{\sim} kQ$ as algebras. Let $\delta: k\Delta\to k\Delta/\langle E\rangle$ be the canonical surjection. Then we have the homomorphisms of algebras$$kQ\lraf{(-)^+}k\Delta\lraf{\delta}k\Delta/\langle E\rangle\lraf{\sim} kQ $$
such that their composition is the identity map of $kQ$. This shows that $(-)^+$ is injective. Applying the map $(-)^+$, we define
$$
\rho^+:=\{\sigma^+\mid \sigma\in\rho\}\quad\mbox{and}\quad \psi':=\rho^+\cup\big(\bigcup_{i,j\in V_0 }(e_ik\Delta e_{\overline{j}}\cup e_{\overline{j}}k\Delta e_i)\big).
$$
We shall show $\langle\psi'\rangle=\langle\psi\rangle$ in $k\Delta$.

In fact, let $\varphi= \bigcup_{i,j\in V_0 }(e_ik\Delta e_{\overline{j}}\cup e_{\overline{j}}k\Delta e_i)\subseteq \psi'$. Clearly, $\langle\varphi\rangle=\langle\psi_1\rangle$. Now, let us consider the image of a path under $(-)^+$.

$(i)$ For $p\in \mathscr{P}(Q)$ of length at least $1$, we have

1) If $p\in\mathscr{P}(Q'),$ then $p^+=p.$

2) If $p\not\in\mathscr{P}(Q'),$ then $p^+=p+\overline{p}+ p'$ with $p'$ in the $k$-space $k\varphi$ generated by elements of $\varphi$.

$(ii)$ For $\sigma\in \rho$, we write $\sigma = \sum_{i=1}^sa_ip_i + \sum_{j= s+1}^na_jp_j$ such that $p_i\in \mathscr{P}(Q')$ for $1\le i\le s$ and $p_{j}\not\in \mathscr{P}(Q')$ for $s+1\le j\le n$. It follows from $(i)$ that
$$ (\ast)\quad \sigma^+=\sum_{i=1}^sa_ip^+_i + \sum_{j= s+1}^na_jp^+_j=\sum_{i=1}^sa_ip_i + \sum_{j= s+1}^na_j(p_j+\overline{p_j}+p'_j)=\sigma + \sum_{j=s+1}^na_j\overline{p_j}+ \sum_{j= s+1}^na_jp'_j$$
If $\sigma\in\psi_2$, then $s=0$ and $\sigma^+=\sigma+\overline{\sigma}+ \sum_{j= 1}^na_jp'_j$ with $\overline{\sigma}\in\psi_3$, and therefore $\sigma^+\in \langle\psi\rangle$. If $\sigma\not\in\psi_2$, that is $s(\sigma), t(\sigma)\in Q_0'$, then $\sigma_+\in\psi_4$ and $\sigma^+=\sigma_+ + \sum_{j= s+1}^na_jp'_j\in \langle\psi\rangle$.  Thus $\langle\psi'\rangle\subseteq\langle\psi\rangle$ in $k\Delta$.

Conversely, pick up $\tau\in \psi$, we show $\tau\in \langle\psi'\rangle.$
If $\tau=\sigma_+\in \psi_4$, then $\tau=\sigma^+ - \sum_{j= s+1}^na_jp'_j\in \langle\psi'\rangle.$ If $\tau = \sigma \in \psi_2$ and $s(\sigma)\in V_0$, then $e_{s(\sigma)}\overline{\sigma}=0$ and therefore $\sigma=e_{s(\sigma)}\sigma=e_{s(\sigma)}\sigma^+-e_{s(\sigma)}\sum_{j=1}^na_jp'_j\in\langle\psi'\rangle.$
If $\tau = \sigma \in \psi_2$ and $t(\sigma)\in V_0$, then $\overline{\sigma}e_{t(\sigma)}=0$ and $\sigma=\sigma e_{t(\sigma)}=\sigma^+e_{t(\sigma)}-\sum_{j=1}^na_jp'_j e_{t(\sigma)}
\in\langle\psi'\rangle.$ If $\tau=\overline{\sigma}\in \psi_3$ with $\sigma\in \psi_2$, then $\overline{\sigma}=\sigma^+-\sigma-\sum_{j= 1}^na_jp'_j$. By what we have just proved, $\sigma\in\langle\psi'\rangle$, and therefore $\overline{\sigma}\in \langle\psi'\rangle$. Thus $\langle\psi\rangle\subseteq \langle\psi'\rangle,$ and therefore $\langle\psi'\rangle=\langle\psi\rangle$ and $k(\Delta,\psi') = k(\Delta,\psi)$.

Since $\varphi\subseteq\langle E \rangle$, it is clear that $\langle\psi'\rangle\subseteq\langle\rho^+\cup E\rangle$.
By the third equality in $(\ast)$ and the fact that $\sum_{j=s+1}^na_j\overline{p_j}$ and $\sum_{j= s+1}^na_jp'_j$
belong to $\langle E\rangle$,  we obtain $\langle\rho^+\cup E\rangle=\langle\rho\cup E\rangle$ in $k\Delta.$
Thus $k\Delta/\langle\rho^+\cup E\rangle=k\Delta/\langle\rho\cup E\rangle\simeq kQ/\langle\rho\rangle= A$ as algebras. Moreover, since $\langle\rho^+\rangle\subseteq\langle\psi'\rangle\subseteq\langle\rho^+\cup E\rangle\subseteq k\Delta,$ the homomorphisms $(-)^+$ and $\delta$ induce algebra homomorphisms $\mu: A\to k\Delta/\langle\psi'\rangle$ and $\overline{\delta}: k\Delta/\langle\psi'\rangle\to k\Delta/\langle\rho^+\cup E\rangle$, respectively.
Now, we identify $k\Delta/\langle\rho^+\cup E\rangle$ with $A$. Then $\mu\overline{\delta}={\rm Id}_A$ and $\mu$ is injective.

$(2)$ We first construct a map $\theta$ by applying Lemma \ref{generator}(3). For simplicity, let
$$R:=R(A,e),\; S:=k(\Delta, \psi),\;x:=\sum_{i\in V_0}e_{\overline{i}}\in S.$$
Then $x^2=x$. By $(1)$, $(e)\mu=e^+=\sum_{j\in V_0}(e_j+e_{\overline{j}})$. Since $e^+ e_{\overline{i}}=e_{\overline{i}}=e_{\overline{i}}e^+$, we have $e^+x=x=xe^+$ and $x\in e^+Se^+$. Recall that $e_jSe_{\overline{i}}=e_{\overline{i}}Se_j=0$ for $i,j\in V_0$, due to the relation set $\psi_1$. Thus, for $s\in S$, we have
$$e^+se^+x=e^+sx=\sum_{j\in V_0}\sum_{i\in V_0}(e_j+e_{\overline{j}})se_{\overline{i}}=\big(\sum_{j\in V_0}e_{\overline{j}}\big)s\big(\sum_{i\in V_0}e_{\overline{i}}\big),$$
$$x e^+se^+=xse^+=\sum_{i\in V_0}\sum_{j\in V_0}e_{\overline{i}}s(e_j+e_{\overline{j}})=\big(\sum_{i\in V_0}e_{\overline{i}}\big)s\big(\sum_{j\in V_0}e_{\overline{j}}\big).$$
This shows $e^+se^+x=xe^+se^+$. Since $\Lambda=eAe$ and $(\Lambda)\mu\subseteq e^+Se^+$, we have $(c)\mu x=xe^+(c)\mu$ for any $c\in\Lambda$.
By Lemma \ref{generator}(3), there is a unique algebra homomorphism $\theta: R\to S$ such that the restriction of $\theta$ to $A$
equals $\mu$ and $(\overline{e})\theta=x$. Let $\overline{e_i}:=e_i\otimes e_i\in R$.
Then $\overline{e_i}=e_i\overline{e}e_i$ and $(\overline{e_i})\theta=e_i^+x e_i^+= e_i^+(\sum_{i\in V_0}e_{\overline{i}})e_i^+ =e_{\overline{i}}$.

Next, we prove that $\theta$ is surjective. It suffices to show that $\Delta_1\subseteq \Img(\theta)$ and $e_t\in \Img(\theta)$ for $t\in\Delta_0$.

In fact, if $t\in Q'_{0}$, then $(e_t)\theta=(e_t)\mu=e_t$; if $t\in V_0$, then
$(\overline{e_t})\theta=e_{\overline{t}}$ and
$(e_t-{\overline{e_t}})\theta=e_t+e_{\overline{t}}-e_{\overline{t}}=e_t$. This implies that $e_t\in\Img(\theta)$ for any $t\in\Delta_0$.
Now, let $\alpha: u\to v\in Q_1$. If $u,v\in Q'_{0}$, then $(\alpha)\theta=\alpha$.
If $u\in V_0$ or $v\in V_0$, then $(\alpha)\theta=(\alpha)\mu=\alpha+\overline{\alpha}$. In case of $u\in V_0$, we get
$$(\overline{e_u}\alpha)\theta=(\overline{e_u})\theta(\alpha)\theta=e_{\overline{u}}(\alpha)\mu=e_{\overline{u}}(\alpha+\overline{\alpha})=\overline{\alpha}\quad\mbox{and}\quad (\alpha-\overline{e_u}\alpha)\theta=\alpha.$$
In case of $v\in V_0$, we have $(\alpha\overline{e_v})\theta=\overline{\alpha}$ and $(\alpha-\alpha\overline{e_v})\theta=\alpha$.
Thus $Q_1\subseteq \Img(\theta)$ and $\overline{Q_1}\setminus Q'_{1}\subseteq \Img(\theta)$.

Finally, we construct an algebra homomorphism $\pi: S\to R$ such that $\theta\pi={\rm Id}_R$,  the identity map of $R$. This means that $\theta$ is injective. Hence  it is bijective.

We define a map $\{e_t\mid t\in \Delta_0\}\cup \Delta_1\to R$ by
$e_i \mapsto e_i-\overline{e_i},\; \;e_{\overline{i}}\mapsto \overline{e_i} \; \mbox{ for } i\in V_0; \quad e_j\mapsto e_j \; \mbox{ for } j\in Q'_{0};$
and for $\alpha\in Q_1$,
$$
i\lraf{\alpha} j\mapsto \left\{
\begin{array}{llll}
\alpha & \hbox{$i, j\in Q'_{0}$},\\
\alpha-\alpha \overline{e_j}, & \hbox{$i\in Q_0, j\in V_0$},\\
\alpha-\overline{e_i}\alpha, & \hbox{$i\in V_0, j\in Q_0;$}
\end{array}\right. \quad
\overline{i}\lraf{\overline{\alpha}} \overline{j}
\mapsto \left\{
\begin{array}{ll}
\alpha \overline{e_j}, & \hbox{$i\in Q_0, j\in V_0$},\\
\overline{e_i}\alpha, & \hbox{$i\in V_0, j\in Q_0,$}
\end{array}\right. \quad
$$
Note that $\overline{e_i}\alpha=e_i\otimes\alpha=\alpha\otimes e_j=\alpha \overline{e_j}$ in $R$ for $i,j\in V_0$.  By Lemma \ref{quiver},  the map can be extended to a unique algebra homomorphism $\gamma: k\Delta\to R$.  Clearly, $\gamma$ preserves the idempotents corresponding to the vertices in $Q'_{0}$ and also the arrows in $Q'_{1}$. Further, if $i\in V_0$,  then $(e_i^+)\gamma=(e_i+e_{\overline i})\gamma=e_i$; if $\alpha\in Q_1\setminus Q'_{1}$,   then $(\alpha^+)\gamma=(\alpha+\overline{\alpha})\gamma=\alpha$. This implies
 $(\sigma^+)\gamma=\sigma$ for any $\sigma\in\rho$. Moreover, by Lemma \ref{Property}(1),
$$(e_ik\Delta e_{\overline{j}})\gamma\subseteq (e_i-\overline{e_i}) R\overline{e_j}\subseteq (e-\overline{e})R\overline{e}=0\;\;\mbox{and}\;\; (e_{\overline j}k\Delta e_i)\varphi\subseteq \overline{e_j} R(e_i-\overline{e_i})\subseteq \overline{e}R(e-\overline{e})=0$$ for any $i,j\in V_0$. Consequently, we have $\langle\psi'\rangle\subseteq \Ker(\gamma)$, and therefore $\gamma$ induces an algebra homomorphism $\pi: S\to R$. Now, let $g:=\theta\pi: R\to R$ and $h:=(-)^+\;\gamma: kQ\to R$.
Since the restriction of $\theta$ to $A$ equals $\mu$, the restriction $g|_A: A\to R$ of $g$ to $A$ is induced from $h.$
As $\gamma$ preserves the idempotents corresponding to the vertices in $Q_0$ and also the arrows in $Q_1$, we see that $g|_A$ has its image in $A$ and  factorizes through  ${\rm Id}_A$.  Since  $(\overline{e_i})g=(e_{\overline{i}})\pi=\overline{e_i}$ for $i\in V_0$ and $\overline{e}=\sum_{i\in V_0}\overline{e_i}$, we have  $(\overline{e})g=\overline{e}$. Thus $g={\rm Id}_R$ by Lemma \ref{generator}(3).  $\square$

Now, let us illustrate the construction $R(A,e)$ by an example.
\begin{Bsp}\label{Exmaple1}
\rm{
Suppose that $A$ is an algebra over a field $k$ presented by the quiver with
relations:
$$
\xymatrix{
1\ar_-{\alpha}@/_0.8pc/[d]\ar^-{\beta}@/^0.8pc/[d]& 4\ar[l]_-{\delta}\ar[d]^-{\sigma} && \\
 2\ar[d]_-{\gamma}\ar[r]^{\tau} & 5\ar@(ur,dr)^-{\eta \;,}\ar[dl]^-{\theta} &&\eta^2=\sigma\eta=\tau\eta=\alpha\gamma=\delta\beta\tau=0, \;\;\beta\gamma=\beta\tau\theta.\\
3&  &}
$$
Let $Q'$ be the full subquiver of $Q$ consisting of the vertex set $\{1,2,3\}$ and let $e=e_4+e_5$. By  Proposition \ref{Quiver algebra}(2), the algebra $R(A, e)$ is isomorphic to the algebra presented by the following quiver with relations:
\vspace{-2cm}
$$\begin{array}{ccc}
\xymatrix{
\overline{4}\ar[r]^-{\overline{\delta}}\ar[d]_-{\overline{\sigma}}& 1\ar_-{\alpha}@/_0.8pc/[d]
\ar^-{\beta}@/^0.8pc/[d]& 4\ar[l]_-{\delta}\ar[d]^-{\sigma}  \\
\overline{5}\ar@(ul,dl)_-{\overline{\eta}}\ar[dr]_-{\overline{\theta}}& 2\ar[d]_-{\gamma}\ar[l]_-{\overline{\tau}}\ar[r]^{\tau} & 5\ar@(ur,dr)^-{\eta \;,}\ar[dl]^-{\theta}\\
&  3&}  & \quad & \begin{array}{l}  \\  \\ \\ \\ \\ \delta\beta\overline{\tau}=\overline{\delta}\beta\tau=\delta\alpha\overline{\tau}=\overline{\delta}\alpha\tau=0,\\
\eta^2=\sigma\eta=\tau\eta=\delta\beta\tau=0,\\
\overline{\eta}^2=\overline{\sigma}\,\overline{\eta}=\overline{\tau}\,\overline{\eta}=\overline{\delta}
\beta\,\overline{\tau}=0,\\ \alpha\gamma=0,\quad \beta\gamma=\beta\tau\theta+\beta\overline{\tau}\,\overline{\theta}.\end{array}\\
\end{array} $$
This quiver is the mirror reflection of the one of $A$ along the full subquiver $Q'$ of $Q$.
}
\end{Bsp}

\section{Mirror-reflective algebras and gendo-symmetric algebras} \label{MSAGSA}

This section is devoted to proofs of all results mentioned in the introduction. We first show that mirror-reflective algebras of gendo-symmetric algebras at any levels are symmetric (see Proposition \ref{Mirror-symmetric}). By iterating this procedure, we construct not only gendo-symmetric algebras of increasing dominant dimensions and higher minimal Auslander-Gorenstein (see Theorem \ref{Main properties}), but also recollements of derived module categories of these algebras (see Theorem \ref{Recollement}). The constructed recollements are then applied to give a new formulation of the Tachikawa's second conjecture for symmetric algebras in terms of stratified dimensions and ratios (see Theorem \ref{Tachikawa's second conjecture}). Consequently, a sufficient condition is given for the conjecture to hold for symmetric algebras (see Corollary \ref{Derived simple}).

Throughout this section, all algebras considered are finite-dimensional algebras over a field $k$.

\subsection{Relations among mirror-reflective, symmetric and gendo-symmetric algebras}\label{SAAGSA}
Let $A$ be an algebra, $e^2=e\in A$ and $\Lambda:=eAe$. Suppose that there is an isomorphism
$\iota: eA\to D(Ae)$ of $\Lambda$-$A$-bimodules. Let $\iota_e:=(e)\iota\in D(Ae)=\Hom_k(Ae,k).$
Then $\iota_e=e\iota_e=\iota_ee$. Moreover, $\iota$ is nothing else than the right multiplication map by $\iota_e$. Define
$\zeta:Ae\otimes_\Lambda eA\to k$ to be the composite of
the maps
$$Ae\otimes_\Lambda eA\lraf{\;id\otimes\iota\;} Ae\otimes_\Lambda D(Ae)\lraf{{\rm ev}}k$$
where ${\rm ev}$  stands for the evaluation map: $ae\otimes f\mapsto (ae)f$ for $a\in A$ and $f\in D(Ae)$. Then
$\zeta$ is given by $(ae\otimes eb)\zeta=(bae)\iota_e=(ebae)\iota_e$ for $a, b\in A$.  Fix an element $\lambda\in Z(\Lambda)$, there are associated two maps
$$\chi:\;\; R(A,e,\lambda)=A\oplus Ae\otimes_\Lambda eA \lra k,\;\; a+\sum_{i=1}^n a_ie\otimes eb_i\;\mapsto\;\sum_{i=1}^n(a_ie\otimes eb_i)\zeta=\sum_{i=1}^n(eb_ia_ie)\iota_e \mbox{ for } a_i, b_i\in A,$$
$$
\gamma:\;\;Ae\otimes_\Lambda eA\lra D(A),\;\; ae\otimes eb\;\mapsto\; [a'\mapsto (eba'ae)\iota_e \, \mbox{ for } a,a', b\in A]
$$
of $k$-spaces. They have the properties.

\begin{Lem} \label{SS}
$(1)$ For any $r_1, r_2\in R(A,e,\lambda)$, $(r_1\ast r_2)\chi=(r_2\ast r_1)\chi$, where $\ast$ denotes the multiplication of $R(A,e,\lambda)$.

$(2)$ The map $\gamma$ is a homomorphism of $A$-$A$-bimodules. It is an isomorphism if and only if the map $(\cdot e): \End_{A\opp}(A)\to\End_{\Lambda\opp}(Ae)$ induced from the right multiplication by $e$ is an isomorphism of algebras.

$(3)$ If $\epsilon:D(A)\to k$ denotes the map sending $f\in D(A)$ to $(1)f$, then $\zeta=\gamma\,\epsilon$.
\end{Lem}

{\it Proof.} $(1)$ It suffices to show
$\big((a_1+ae\otimes eb)\ast(a_2+a'e\otimes eb')\big)\chi=\big((a_2+a'e\otimes eb')\ast(a_1+ae\otimes eb)\big)\chi$
for any $a, a', b,b', a_1,a_2\in A$. However, this follows from
$\big(a'(ae\otimes eb)\big)\zeta=\big((ae\otimes eb)a'\big)\zeta$ and $
\big((ae\otimes eb)\otimes(a'e\otimes eb')\big)\omega_\lambda\zeta=\big((a'e\otimes eb')\otimes(ae\otimes eb)\big)\omega_\lambda\zeta,
$ by the definitions of $\zeta$ and $\omega_\lambda$ (see Section \ref{sect3.1} for definition).

$(2)$ Note that there is a canonical isomorphism
$\varphi:\;\; Ae\otimes_\Lambda D(Ae)\ra D(\End_{\Lambda\opp}(Ae)),\;\; ae\otimes f\mapsto [g\mapsto (ae)gf]$
for $a\in A$, $f\in D(Ae)$ and $g\in \End_{\Lambda\opp}(Ae)$. Let $\vartheta: A\to\End_{A\opp}(A)$ be the isomorphism which sends $a$ to $(a\cdot)$.
Then the composition of the maps
$$
Ae\otimes_\Lambda eA\lraf{\;Ae\otimes\iota\;} Ae\otimes_\Lambda D(Ae)\lraf{\varphi}D(\End_{\Lambda\opp}(Ae))\lraf{D(\cdot e)} D(\End_{A\opp}(A))\lraf{D(\vartheta)}D(A)
$$
coincides with $\gamma$. Clearly, all the maps above are homomorphisms of $A$-$A$-bimodules.
Thus $\gamma$ is a homomorphism of $A$-$A$-bimodules. Since $D:k\modcat\to k\modcat$ is a duality, $\gamma$ is an isomorphism if and only if the map $(\cdot e)$ in $(2)$ is an isomorphism of algebras.

$(3)$ This follows from $(ae\otimes eb)\zeta=(ebae)\iota_e$ for $a,b\in A$. $\square$

\medskip
\textbf{From now on}, let $(A, e)$ be a \emph{gendo-symmetric} algebra. Recall that  $\add(Ae)$ coincides with the full subcategory of $A\modcat$
consisting of projective-injective $A$-modules. If $e'$ is another idempotent of $A$ such that
$\add(Ae)=\add(Ae')$, then the mirror-reflective algebras $R(A,e)$ and $R(A,e')$ are isomorphic as algebras by Lemma \ref{generator}(1).
So, for simplicity, we write $R(A)$ for $R(A,e).$

In the following, we describe $R(A)$ as deformation of a trivial extension.
Let $\Lambda:=eAe$ and  $\iota: eA\to D(Ae)$ be an isomorphism of $\Lambda$-$A$-bimodules (see Lemma \ref{GS}(2)).
Then $\Lambda$ is symmetric and $eA$ is a generator over $\Lambda$.
Moreover, there are algebra isomorphisms $A\simeq\End_\Lambda(eA)$ and $A\opp\simeq\End_{\Lambda\opp}(Ae)$.
By Lemma \ref{SS}(2), there is an isomorphism of $A$-$A$-bimodules:
$\gamma:\;\;Ae\otimes_\Lambda eA\lraf{\simeq} D(A).$ Since $A\simeq \End_\Lambda(eA)$ and $eA$ is a generator over $\Lambda$,
the functor $e(-)e: A^{e}\Modcat\to\Lambda^e\Modcat$ between the categories of bimodules induces an algebra isomorphism
$Z(A)\to Z(\Lambda)$. So, for $\lambda\in Z(\Lambda)$, there exists a unique element $\lambda'\in Z(A)$ such that $e\lambda'e=\lambda$.
Define  $\overline{\omega_e}:=(\gamma\otimes\gamma)^{-1}\omega_e\gamma:\;\; D(A)\otimes_AD(A)\lraf{\simeq} D(A)$ and $F =Ae\otimes_\Lambda-\otimes_\Lambda eA: \Lambda^{e}\Modcat\to A^{e}\Modcat.$
We obtain the commutative  diagram
$$
\xymatrix{(Ae\otimes_\Lambda eA)\otimes_A(Ae\otimes_\Lambda eA)\ar[r]^-{\omega_e}_-{\simeq}\ar[d]_-{\gamma\otimes\gamma}&Ae\otimes_\Lambda eA\ar[r]^-{F(\cdot\lambda)}\ar[d]^-{\gamma}& Ae\otimes_\Lambda eA\ar[d]^-{\gamma}\\
D(A)\otimes_AD(A)\ar[r]^-{\overline{\omega_e}}& D(A)\ar[r]^-{(\cdot\lambda')}& D(A).}
$$
Define $\overline{\omega_\lambda}:=\overline{\omega_e}(\cdot\lambda'): \;\; D(A)\otimes_AD(A)\lra D(A).$
Now, we extend  $\overline{\omega_\lambda}$ to a multiplication on the direct sum $A\oplus D(A)$ by setting
$$(A\oplus D(A))\times (A\oplus D(A))\lra A\oplus D(A),\;\;\big((a, f), (b,g)\big)\mapsto \big(ab, ag+fb+(f\otimes g)\overline{\omega_\lambda}\big)$$
for $a,b\in A$ and $f,g\in D(A)$. Denote by $A\ltimes_\lambda D(A)$ the abelian group $A\oplus D(A)$ with the above-defined multiplication.
By Lemma \ref{rc}(1), $A\ltimes_\lambda D(A)$ is an algebra with an algebra isomorphism
$$\overline{\gamma}:=\left(\begin{array}{lc}{\rm Id}_A &0 \\
0 & \gamma\end{array}\right):\;\; R(A,e,\lambda)\lraf{\simeq} A\ltimes_\lambda D(A).$$
Compared with the trivial extension $A\ltimes D(A)$, the following result, suggested by Kunio Yamagata, shows that
$A\ltimes_\lambda D(A)$ is also a symmetric algebra for any $\lambda$.

\begin{Prop} \label{Mirror-symmetric}
If $(A, e)$ is a gendo-symmetric algebra, then $R(A,e,\lambda)$ is symmetric for $\lambda\in Z(\Lambda)$.
\end{Prop}

{\it Proof.}  Let $R:=R(A,e,\lambda)$. Applying $\chi:R\to k$, we define a bilinear form $\widetilde{\chi}: R\times R\to k$, $(r_1,r_2)\mapsto (r_1\ast r_2)\chi$ for $r_1, r_2\in R$. By Lemma \ref{SS}(1), $\widetilde{\chi}$ is symmetric. To show that $R$ is a symmetric algebra, it suffices to show that $\widetilde{\chi}$ is non-degenerate.

Let $T:=A\ltimes_\lambda D(A)$ and $\psi:=\overline{\gamma}^{-1}\chi:T\to k$. Since $\overline{\gamma}:R\to T$ is an algebra isomorphism,
$\psi$  induces a symmetric bilinear form $\widetilde{\psi}:T\times T\to k$, $(t_1, t_1)\in T\times T\mapsto (t_1t_2)\psi$.
Clearly, $\widetilde{\chi}$ is non-degenerate if and only if so is $\widetilde{\psi}$. Further, by Lemma \ref{SS}(3), $\psi$
is given by $(a, f)\mapsto (1)f$ for $a\in A$ and $f\in D(A)$. This implies that
$\big((a, f), (b,g)\big)\widetilde{\psi}=(a)g+(b)f+(1)(f\otimes g)\overline{\omega_\lambda}$ for $b\in A$ and $g\in D(A)$.
Now, we show that $\widetilde{\psi}$ is non-degenerate.

Let $(a,f)\neq 0$. Then $a\neq 0$ or $f\neq 0$.  If $f\neq 0$, then there is an element $b\in A$ such that $(b)f\neq 0$, and therefore $\big((a, f), (b,0)\big)\widetilde{\psi}=(b)f\neq 0$. If $f=0$ and $a\neq 0$, then the canonical isomorphism $A\simeq DD(A)$ implies that there is an element $g\in D(A)$
such that $(a)g\neq 0$. In this case, $\big((a, 0), (0,g)\big)\widetilde{\psi}=(a)g\neq 0$. Thus $\widetilde{\psi}$ is non-degenerate.
$\square$

\medskip
Compared with $R(A)$, the algebra $S(A,e)$ depends on the choice of $e$, that is, if $f=f^2\in A$ such that $(A,f)$ is gendo-symmetric, then $S(A,e)$ and $S(A,f)$ do not have to be isomorphic in general.
The following result collects basic homological properties of $S(A,e)$.

\begin{Prop}\label{reduced MS}
Let $S:=S(A,e)$ and $B_0:=(1-e)A(1-e)$. Then

$(1)$ $S$ is a symmetric algebra.

$(2)$ $B_0$ can be regarded as a $S$-module and contains no nonzero projective direct summands.

$(3)$ If $\add({_A}Ae)\cap\add({_A}A(1-e))=0$, then $\#(S)=\#(A)$. For instance, if $B_0$ is indecomposable as an algebra, then so is $S$.
\end{Prop}

{\it Proof.} $(1)$ Let $R:=R(A)$, $\bar{e}:=e\otimes e\in R$ and $e_0:=(1-e)+\overline{e}\in R$. Since $R$ is symmetric by Proposition \ref{Mirror-symmetric}(1) and $S=e_0Re_0$ by Lemma \ref{Property}(1), $S$ is symmetric.

$(2)$ Since $\pi_1$ induces a surjective algebra homomorphism $\pi_1': S\to B_0$ such that $S/S\overline{e}S\simeq B_0$ (see Lemma \ref{Property} for notation), $B_0$ can be regarded as an $S$-module. Assume that the $S$-module $B_0$ contains an indecomposable projective direct summand $X$.
Then there is a primitive idempotent $f\in A$ such that $1-e=f+f'$ with $f$ and $f'$ orthogonal idempotents in $A$, and $X\simeq Sf$ as $S$-modules.
Clearly, $S\overline{e}Sf=0$, $(f)\pi_2'=f$, $(1-e)\pi_2'=1-e$ and $(S\overline{e}Sf)\pi_2'=AeAf$.
Consequently,  $\Hom_A(Ae,Af)\simeq eAf=0$, and therefore $\Hom_A(Af, Ae)\simeq D\Hom_A(Ae,Af)=0$. By Lemma \ref{GS}(2), $Af$ can be embedded into $(Ae)^n$ for some $n\geq 1$. This implies $Af=0$, a contradiction.

$(3)$ Since $\bar{e}S\bar{e}\simeq eAe$ by Lemma \ref{ring homo}(2),  it follows from $(2)$ that $\#S(A)=\#(eAe)+\#(B_0)$. Due to  $\add(Ae)\cap\add(A(1-e))=0$, we have $\#(A)=\#(eAe)+\#(B_0)$ and $\#S(A)=\#(A)$. The second assertion in $(3)$ follows from Proposition \ref{Indecomposable}(2). $\square$

\subsection{Mirror-reflective algebras and Auslander-Gorenstein algebras}

In the subsection, we construct new gendo-symmetric algebras from minimal Auslander-Gorenstein algebras and then present a proof of Theorem \ref{Main properties}. This is based on study of mirror-reflective algebras.

Recall from Lemma \ref{Property} that we have an algebra automorphism $\phi:R(A)\to R(A)$ and two surjective algebra homomorphisms $\pi_1, \pi_2: R(A)\to A$ such that $\pi_2=\phi\pi_1$. Thus we  regard $A$-modules as $R(A)$-modules via $\pi_1$ in the following discussion. It turns out that $A\modcat$ is a Serre subcategory of $R(A)\modcat$, that is, it is closed under direct summands, submodules, quotients and extensions in $R(A)\Modcat$. Let
$$\phi_*: R(A)\modcat\lra R(A)\modcat \quad \mbox{and} \quad (\pi_2)_*:A\modcat\lra R(A)\modcat $$
be the restriction functors induced by $\phi$ and $\pi_2$, respectively.  Then $\phi_*$ is an auto-equivalence and $\phi_*(X)=(\pi_2)_*(X)$ for each $A$-module $X$.

\begin{Lem}\label{Periodic}
Suppose that $\Lambda$ is a symmetric algebra and $N$ is a basic $\Lambda$-module without nonzero projective direct summands. Let $A:=\End_\Lambda(\Lambda\oplus N)$, $e$ an idempotent of $A$ corresponding to the direct summand $\Lambda$ of $\Lambda\oplus N$, and $R:=R(A, e)$. If ${_\Lambda}N$ is $m$-rigid for a natural number $m$, then the following hold.

$(1)$ The $R$-module $A(1-e)$ is $(m+2)$-rigid and there are isomorphisms of $R$-modules:
$$\Omega_R^{m+3}\big(A(1-e)\big)\simeq\Omega_R^{m+2}\big(\phi_*(Ae\otimes_\Lambda N)\big)\simeq \phi_*\big(\Hom_\Lambda(eA, \Omega_\Lambda^{m+2}(N))\big).$$

$(2)$ If $\, \Omega_\Lambda^{m+2}(N)\simeq N$, then $\Omega_R^{m+3}(A(1-e))\simeq \phi_*(A(1-e))$ and the $R$-module $A(1-e)$ is $(2m+4)$-rigid. In this case, $\Omega_R^{2m+6}(A(1-e))\simeq A(1-e)$.
\end{Lem}

{\it Proof.} $(1)$ By the proof of Proposition \ref{Idem-rec}(2), $\pi_2$ induces an isomorphism
$\Tor_i^{\overline{e}R\overline{e}}(R\overline{e}, \overline{e}R)\simeq \Tor_i^\Lambda(Ae, eA)$
for all $i\geq 1$. Since $\Lambda$ is symmetric and $D(Ae)\simeq eA$ by Lemma \ref{GS}(2), we have
$$
D\Tor_i^\Lambda(Ae, eA)\simeq\Ext_\Lambda^i(eA, D(Ae))\simeq\Ext_\Lambda^i(eA, eA)=\Ext_\Lambda^i(\Lambda\oplus N, \Lambda\oplus N)\simeq\Ext_\Lambda^i(N,N).
$$
As ${_\Lambda}N$ is $m$-rigid, there holds $\Tor_i^{\overline{e}R\overline{e}}(R\overline{e}, \overline{e}R)=0$ for $1\leq i\leq m$.
By Proposition \ref{Idem-rec}(1), $I:=R\overline{e}R$ is $2$-idempotent. Therefore
$I$ is $(m+2)$-idempotent by Lemma \ref{strong idempotent}(1). Further, it follows from Lemma \ref{strong idempotent}(2) that ${_R}R/I$ is $(m+2)$-rigid.
Since $R/I\simeq A$ as $R$-modules,  $_RA$ is $(m+2)$-rigid. Note that ${_R}A\simeq R(e-\overline{e})\oplus A(1-e)$ by Lemma \ref{ring homo}(2).
As $R$ is symmetric by Proposition \ref{Mirror-symmetric}, we see that $R(e-\overline{e})$ is projective-injective. Consequently, ${_R}A(1-e)$ is $(m+2)$-rigid.

The proof of Proposition \ref{Idem-rec}(1) implies $I\simeq R\overline{e}\otimes_{\overline{e}R\overline{e}}\overline{e}R$ as $R$-$R$-bimodules. By
Lemma \ref{ring homo}(2), $\pi_2$ restricts to an algebra isomorphism $\overline{e}R\overline{e}\to\Lambda$ and also an isomorphism $R\overline{e}\to Ae$ of abelian groups.  Via the algebra isomorphism, we can regard $R\overline{e}$ as an $R$-$\Lambda$-bimodule.
Then $R\overline{e}\simeq(\pi_2)_*(Ae)=\phi_*(Ae)$ as $R$-$\Lambda$-bimodules. This gives a natural isomorphism
$R\overline{e}\otimes_\Lambda-\lraf{\simeq} \phi_*(Ae)\otimes_\Lambda-$ of functors  from $\prj{\Lambda}$ to $\prj{R}$.
Since $N$ has no nonzero projective direct summands, $\add({_A}Ae)\cap\add({_A}A(1-e))=0$.
From $A\otimes_RR\overline{e}\simeq Ae\simeq R\overline{e}$ and $A\otimes_RR(1-e)\simeq A(1-e)$,  we obtain $\add(R\overline{e})\cap\add(R(1-e))=0$. Since $I(1-e)$ is isomorphic to $R\overline{e}\otimes_{\overline{e}R\overline{e}}\overline{e}R(1-e)$ which is a quotient module of $(R\overline{e})^n$ for some $n$,  we deduce that $I(1-e)$ does not contain nonzero direct summands in $\add(R(1-e))$. Thus
the surjection ${_R}R(1-e)\to A(1-e)$ induced by $\pi_1$ is a projective cover of the $R$-module $A(1-e)$, and therefore $\Omega_R(A(1-e))=I(1-e)$.
Since $\pi_2$ induces an isomorphism $\overline{e}R\to eA$ and sends $1-e$ to $1-e$ by Lemma \ref{ring homo}(2),
we have $\overline{e}R(1-e)\simeq eA(1-e)$ and
$$\Omega_R(A(1-e))\simeq R\overline{e}\otimes_{\overline{e}R\overline{e}}eA(1-e)\simeq R\overline{e}\otimes_\Lambda eA(1-e)\simeq \phi_*(Ae)\otimes_\Lambda N=\phi_*(Ae\otimes_\Lambda N).
$$
Let $\cdots\to Q_{m+1}\lraf{\partial} Q_m\to\cdots\to Q_1\to Q_0\to N\to 0$ be a minimal projective resolution of $_\Lambda N$. Then it follows from $eA(1-e)=N$ and $\Tor_i^\Lambda(Ae, N)\simeq D\Ext_\Lambda^i(N,N)=0$ for $1\leq i\leq m$ that the sequence $$Ae\otimes_\Lambda Q_{m+1}\lraf{Ae\otimes\partial} Ae\otimes_\Lambda Q_m\lra\cdots\lra Ae\otimes_\Lambda Q_1\lra Ae\otimes_\Lambda Q_0\lra Ae\otimes_\Lambda N\lra 0$$ is exact. As the composition of ${_A}Ae\otimes_\Lambda-$ with $(e\cdot)$ is isomorphic to the identity functor of $\Lambda\modcat$, we have $\Omega_A^{m+2}(Ae\otimes_\Lambda N)\simeq \Ker(Ae\otimes\partial)$. Note that $Ae\otimes_\Lambda-\simeq\Hom_\Lambda(eA,-): \prj{\Lambda}\lraf{\simeq}\add({_A}Ae)$ since
$Ae=\Hom_\Lambda(\Lambda\oplus N, \Lambda)$.
This shows  $\Ker(Ae\otimes\partial)\simeq \Hom_\Lambda(eA,\Ker(\partial))=\Hom_\Lambda(eA,\Omega_\Lambda^{m+2}(N))$, and therefore
$$
\Omega_R^{m+3}\big(A(1-e)\big)\simeq\Omega_R^{m+2}(\phi_*(Ae\otimes_\Lambda N))\simeq \phi_*(\Omega_R^{m+2}(Ae\otimes_\Lambda N))\simeq\phi_*(\Hom_\Lambda(eA,\Omega_\Lambda^{m+2}(N))).
$$

$(2)$ Let $X:=A(1-e)$. Suppose $\Omega_\Lambda^{m+2}(N)\simeq N$. Then
$\Omega_R^{m+3}\big(X\big)\simeq\phi_*(\Hom_\Lambda(eA,eX)).$
Since the functor $(e\cdot): A\modcat\to\Lambda\modcat$ induces an algebra isomorphism $\End_A(A)\simeq\End_\Lambda(eA)$, we have $X\simeq\Hom_A(A,X)\simeq\Hom_\Lambda(eA,eX)$. It follows that $\Omega_R^{m+3}\big(X\big)\simeq\phi_*(X)$.
Note that $\phi$ is an algebra isomorphism with $\phi^2={\rm Id}_R$ by Lemma \ref{Property}(3). Since $\Omega_R$ commutes with $\phi_*$, we obtain
$\Omega_R^{2m+6}\big(X\big)\simeq X$. Now, it remains to show that ${_R}X$ is $(2m+4)$-rigid.

Since $R$ is symmetric,  the stable module category $\stmc{R}$ of $R$
is a triangulated category with the shift functor $[1]=\Omega^{-}_R: \stmc{R}\to \stmc{R}$, where  $\Omega^{-}_R$ is the cosyzygy functor on $\stmc{R}$.
Clearly, $\Ext_R^n(X_1, X_2)\simeq\StHom_R(X_1, X_2[n])$ for all $n\geq 1$ and $X_1,X_2\in R\modcat$,
where $\StHom_R(X,Y)$ denotes the morphism set from $X$ to $Y$ in $R\stmc$. Since the Auslander-Reiten (AR) translation on $R\modcat$ coincides with $\Omega^2_R$, it follows from the AR-formula that there is a natural isomorphism $D\StHom_R(X_1, X_2)\simeq\StHom_R(X_2, X_1[-1])$.
Consequently, for each $i\in\mathbb{N}$, there are isomorphisms
$$
\Ext_R^{m+3+i}(X,X)\simeq\StHom_R(\Omega_R^{m+3}(X), X[i])\simeq\StHom_R(\phi_*(X), X[i])\simeq D\StHom_R(X[i],\phi_*(X)[-1]).
$$
Recall that $\phi$ is an algebra isomorphism with $\phi^2={\rm Id}_R$ by Lemma \ref{Property}(3). Then
$$
\StHom_R(X[i],\phi_*(X)[-1])\simeq \StHom_R(\phi_*(X)[i], X[-1])\simeq \StHom_R(\Omega_R^{m+3}(X), X[-1-i])\simeq \Ext_R^{m+2-i}(X,X)
$$
for $0\le i\leq m+1$. This implies $\Ext_R^{m+3+i}(X,X)\simeq D\Ext_R^{m+2-i}(X,X)$ for $0\leq i\leq m+1$. Since $X$ is $(m+2)$-rigid by $(1)$, it is
actually $(2m+4)$-rigid. $\square$

\begin{Prop}\label{OSMO}
Suppose that $\Lambda$ is a symmetric algebra and $N$ is a basic $\Lambda$-module without nonzero projective direct summands. Let  $A :=\End_\Lambda(\Lambda\oplus N)$, $e$ an idempotent of $A$ corresponding to the direct summand $\Lambda$ of $\Lambda\oplus N$, and $R:=R(A, e)$.

$(1)$ If ${_\Lambda}\Lambda\oplus N$ is $m$-rigid, then ${_R}R\oplus A(1-e)$ is $(m+2)$-rigid.

$(2)$ If ${_\Lambda}\Lambda\oplus N$  is $m$-ortho-symmetric, then ${_R}R\oplus A(1-e)$ is $(2m+4)$-ortho-symmetric.

$(3)$ If ${_\Lambda}\Lambda\oplus N$  is maximal $m$-orthogonal, then ${_R}R\oplus A(1-e)$ is maximal $(2m+4)$-orthogonal.
 \end{Prop}

{\it Proof.}  $(1)$ Since $R$ is a symmetric algebra by Proposition \ref{Mirror-symmetric}, (1) follows from Lemma \ref{Periodic}(1).

$(2)$ By assumption, ${_\Lambda}N$ is basic and contains no nonzero projective direct summands. This implies that ${_A}A(1-e)$ is basic and contains no nonzero projective-injective direct summands. We claim that ${_R}A(1-e)$ contains no nonzero projective direct summands.
In fact, by the proof of Lemma \ref{Periodic}(1), ${_R}R(1-e)$ is a projective cover of ${_R}A(1-e)$.
If ${_R}A(1-e)$ contains an indecomposable projective direct summand $Y$, then $Y$ is a direct summand of $R(1-e)$.
Since $R$ is symmetric, ${_R}Y$ must be projective-injective. However, since $A\modcat\subseteq R\modcat$ is a Serre subcategory,
${_A}Y$ is also a nonzero projective-injective direct summand of ${_A}A(1-e)$. This is a contradiction and shows that the above claim holds.
Now $(2)$ follows from Lemmas \ref{Periodic} and \ref{OSS}.

$(3)$ Recall that maximal orthogonal modules over an algebra $B$
are exactly ortho-symmetric $B$-modules such that their endomorphism algebras have finite global dimension.
Let $A_1:=\End_R(R\oplus A(1-e))$. By $(2)$, to show $(3)$, it suffices to show that $\gd(A_1)<\infty$ if $\gd(A)<\infty$

Let $B_1:=\End_R(R\oplus A)$. Since ${_R}A\simeq R(e-\overline{e})\oplus A(1-e)$ by the proof of Lemma \ref{Periodic}(1), $A_1$ and $B_1$ are Morita equivalent, and therefore $\gd(A_1)=\gd(B_1)$.  Since $eA_A$ is faithful, it follows from Proposition \ref{Faithful}(3) that if $\gd(A)<\infty$ then $\gd(B_1)=\gd(B_1\opp)<\infty$. Hence $\gd(A_1)<\infty$. $\square$

\medskip
{\bf Proof of Theorem \ref{Main properties}.}  The statement $(1)$ follows from Proposition \ref{Mirror-symmetric}.
Let $R:=R(A)$ and $S:=S(A,e)$. Then $R$ and $S$ are symmetric by $(1)$ and Proposition \ref{reduced MS}(1). Let $A_2:=\mathcal{A}(A,e)$ and $B_2:=\mathcal{B}(A,e)$.
Then $A_2$ and $B_2$ are gendo-symmetric.

Next, we show that $(2)$ and $(3)$ hold for $A_2$. In fact, since $A$ is gendo-symmetric, we can identify $A$ with $\End_\Lambda(\Lambda\oplus X)$, where $\Lambda:=eAe$ is symmetric and $X=eA(1-e)$.
As global, dominant and injective dimensions are invariant under Morita equivalences, the classes of minimal Auslander-Gorenstein algebras and of higher Auslander algebras are closed under Morita equivalences. Moreover, for a self-injective algebra $\Gamma$ and $M\in \Gamma\modcat$, it follows from \cite[Lemma 3]{Muller} that $\dm(\End_\Gamma(\Gamma\oplus M))$ equals the maximal natural number $n\geq 2$ or $\infty$ such that $M$ is $(n-2)$-rigid. So, for a basic module $X$ that has no nonzero projective direct summands, the inequality $\dm(A_2)\geq \dm(A)+2$ and the statement $(3)$ follow immediately from Proposition \ref{OSMO}. Further, for an arbitrary module $X$, the consideration can be reduced by a series of Morita equivalences, as shown below.

We take a direct summand $N$ of $X$ such that $N$ is basic, has no nonzero projective direct summands and satisfies $\add(\Lambda\oplus N)=\add(\Lambda\oplus X)$. Let $B:=\End_\Lambda(\Lambda\oplus N)$ and $f^2=f\in A$ correspond to the direct summand $\Lambda\oplus N$ of $\Lambda\oplus X$. Then ${_A}Af$ is a progenerator (that is, a projective generator), and therefore $B=fAf$ is Morita equivalent to $A$. Since $ef=e=fe$,
we have $R(B)=fAf\oplus fAe\otimes_\Lambda eAf=fRf$. Due to $R\otimes_AAf\simeq Rf$, the module ${_R}Rf$ is a progenerator. Thus $R$ and $R(B)$ are Morita equivalent.
Now, let $B_2:=\End_{R(B)}(R(B)\oplus B(f-e))$. If $A$ is $n$-minimal Auslander-Gorenstein (respectively, $n$-Auslander), then so is $B$, and therefore, so is $B_2$
by the above-proved case. Next, we shall show that $A_2$ and $B_2$ are Morita equivalent.
Recall that the restriction of $\pi_1$ to $A$ is the identity map of $A$.
This implies $A\otimes_RRf=Af$ as $R$-modules, and therefore $\add(_RA)=\add(_RAf)$. Let $A_2':=\End_R(Rf\oplus A(1-e)f)=\End_R(Rf\oplus A(f-e))$. Then $A_2$ and $A_2'$ are Morita equivalent. Since the functor $(f\cdot): R\modcat\to R(B)\modcat$ is an equivalence and $f(Rf\oplus A(f-e))=R(B)\oplus B(f-e)$, there is an algebra isomorphism $A_2'\simeq B_2$. Hence  $A_2$ and $B_2$ are Morita equivalent.  Thus $(2)$ and $(3)$ hold true for $A_2$.

It remains to show $\dm(B_2)\geq \dm(A)+2$.
Up to Morita equivalence, we assume  $A=\End_\Lambda(\Lambda\oplus N)$. If ${_\Lambda}\Lambda\oplus N$ is $m$-rigid for some $m\in\mathbb{N}$, then it follows from the first part of the proof of Lemma \ref{Periodic}(1) that $I$ is an $(m+2)$-idempotent ideal of $R$. Let  $e_0:=(1-e)+\overline{e}\in R$. By Lemma \ref{Property},
we have $\overline{e}e_0=\overline{e}=e_0\overline{e}$, $I:=R\overline{e}R=S\overline{e}S$ and $S/I\simeq (1-e)A(1-e)$ as algebras. Thanks to Corollary \ref{HID}(1), $I$ is an $(m+2)$-idempotent ideal of $S$.
Further, by Lemma \ref{strong idempotent}(2), ${_S}S/I$ is $(m+2)$-rigid, and therefore ${_S}S\oplus S/I$ is $(m+2)$-rigid since $S$ is symmetric by Proposition \ref{reduced MS}(1). Thus  $\dm(B_2)\geq \dm(A)+2$, due to  \cite[Lemma 3]{Muller}.
$\square$

\subsection{Recollements of mirror-reflective algebras and Tachikawa's second conjecture}\label{IMSAT}

In this subsection, we study the iterated process of constructing (reduced) mirror symmetric algebras from gendo-symmetric algebras and prove Theorems \ref{Recollement} and \ref{Tachikawa's second conjecture}.

Throughout this section, let $(A,e)$ be a gendo-symmetric algebra. We define inductively for $n\ge 1$
$$A_1=B_1:=A,\quad R_1:=R(A_1,e_1), \quad \quad S_1:=S(A_1, f_1),$$
$$ \quad\quad\;\; A_{n+1}:=\End_{R_n}\big(R_n\oplus A_n(1_{A_n}-e_n)\big),\quad R_{n+1}:=R(A_{n+1}, e_{n+1}),$$
$$ \quad \quad\quad\quad\quad\quad B_{n+1}:=\End_{S_n}\big(S_n\oplus (1_{B_n}-f_n)B_n(1_{B_n}-f_n)\big),\quad S_{n+1}:=S(B_{n+1}, f_{n+1}),$$
where $e_1=f_1:=e$, and for $n\geq 1$, $e_{n+1}\in A_{n+1}$ is the idempotent corresponding to the direct summand $R_n$ of the $R_n$-module $R_n\oplus A_n(1_{A_n}-e_n)$, and  $f_{n+1}\in B_{n+1}$ is the idempotent corresponding to the direct summand $S_n$ of the $S_n$-module $S_n\oplus (1_{B_n}-f_n)B_n(1_{B_n}-f_n)$.
In other words, $$A_{n+1}=\mathcal{A}(A_n, e_n),\quad B_{n+1}=\mathcal{B}(B_n, f_n)\, \mbox{ for } n\ge 1.$$
(see Introduction for notation). For convenience, we set $R_0=S_0:=eAe$ and $B_0:=(1-e)A(1-e).$

\begin{Def}
For $n\geq 1$, the algebras $R_n$, $S_n$, $A_n$ and $B_n$ are called the $n$-th mirror-reflective, reduced mirror-reflective, gendo-symmetric and reduced gendo-symmetric algebras of $(A,e)$, respectively.
\end{Def}

By Propositions \ref{Mirror-symmetric} and \ref{reduced MS}(1), the algebras $R_n$ and $S_n$ are symmetric. Thus $A_n$ and $B_n$ are gendo-symmetric. They are characterized in terms of Morita context algebras in Section \ref{SRII}. Moreover, it follows from Theorem \ref{Main properties}(2) that $\dm(A_{n+1})\geq \dm(A_n)+2$ and $\dm(B_{n+1})\geq \dm(B_n)+2$. Thus
$\min\{\dm(A_n),\dm(B_n)\}\geq \dm(A)+2(n-1)\geq 2n.$

\begin{Lem} \label{GMCA}
$(1)$ Let $I_n:=R_n\overline{e}_n R_n$ and $J_n:=R_n(e_n-\overline{e}_n)R_n$ with $\overline{e}_n=e_n\otimes e_n\in R_n$ for $n\geq 1$. Then $A_{n+1}$ is
derived equivalent and stably equivalent of Morita type to the Morita context algebra $M_l(R_n,I_n,J_n)$.

$(2)$ Let $K_n:=S_n\overline{f}_n S_n$ and $L_n:=S_n\cap \big(R(B_n)(f_n-\overline{f}_n)R(B_n)\big)$ for $n\geq 1$.
Then $B_{n+1}$ is derived equivalent and stably equivalent of Morita type to the Morita context algebra $M_l(S_n,K_n, L_n)$.
\end{Lem}

{\it Proof.} $(1)$ Recall that there is a surjective algebra homomorphism  $\pi_{1,n}: R_n\to A_n$ with $\Ker(\pi_{1,n})=I_n$ which induces an isomorphism $R_n(e_n-\overline{e}_n)\simeq A_ne_n$ of $R_n$-modules. Thus $I_n\simeq \Omega_{R_n}(A_n)\oplus Q_n$ with $Q_n$ a projective $R_n$-module, and $A_ne_n$ is a projective $R_n$-module. Hence $A_{n+1}$ is Morita equivalent to $A_{n+1}':=\End_{R_n}(R_n\oplus A_n)$.
Let $C_{n+1}:=\End_{R_n}(R_n\oplus I_n)$. Since $R_n$ is symmetric, it follows from \cite[Corollary 1.2]{HX13} and \cite[Theorem 1.1]{HX10} that $A_{n+1}'$ and $C_{n+1}$ are both derived equivalent and stably equivalent of Morita type. Consequently, $A_{n+1}$ and $C_{n+1}$ are both derived equivalent and stably equivalent of Morita type. It remains to show $C_{n+1}\simeq M_l(R_n,I_n,J_n)$ as algebras.

In fact, since $I_n^2=I_n$, the inclusion $\lambda_n:I_n\hookrightarrow R_n$ induces $\End_{R_n}(I_n)\simeq\Hom_{R_n}(I_n, R_n)$.
As $R_n$ is symmetric and $J_n={\rm Ann}_{R_n\opp}(I_n)$ by Lemma \ref{ring homo}(1), we get $R_n/J_n\simeq\End_{R_n}(I_n)$ as algebras via the restriction of $\lambda_n$. This yields a series of isomorphisms
$$C_{n+1}\simeq\left(\begin{array}{lc} R_n & I_n\\
\Hom_{R_n}(I_n, R_n) & \End_{R_n}(I_n)\end{array}\right)\simeq \left(\begin{array}{lc} R_n & I_n\\
\End_{R_n}(I_n)& \End_{R_n}(I_n)\end{array}\right)\simeq \left(\begin{array}{lc} R_n & I_n\\
R_n/J_n& R_n/J_n\end{array}\right),$$
of which the composition is  an isomorphism from $C_{n+1}$ to $M_l(R_n,I_n,J_n)$ of algebras. This shows $(1)$.

$(2)$ By Lemma \ref{Property}(4), $K_n=R(B_n)\overline{f_n}R(B_n)$ and $S_n/K_n\simeq (1_{B_n}-f_n)B_n(1_{B_n}-f_n)$. By the proof of Proposition \ref{Faithful}(1), ${\rm Ann}_{S_n\opp}(K_n)=L_n$.  Similarly, since $S_n$ is symmetric, we can show that $B_{n+1}$ and $\End_{S_n}(S_n\oplus K_n)$ are both derived equivalent and stably equivalent of Morita type, and that $\End_{S_n}(S_n\oplus K_n)$ is isomorphic to $M_l(S_n,K_n, L_n)$ as algebras.  $\square$

\begin {Rem}\label{Unbounded}
By the proof of Lemma $\ref{GMCA}$,  $B_{n+1}$ and $\End_{S_n}(S_n\oplus \, S_n/K_n)$ are isomorphic, while  $A_{n+1}$ and $\End_{R_n}(R_n\oplus A_n)$ are Morita equivalent. It follows from Proposition \ref{Faithful}(1) that there are recollements of derived module categories
$\big(\D{A_n},\D{A_{n+1}}, \D{A_n}\big)$ and $\big(\D{B_n}, \D{B_{n+1}}, \D{B_0}\big)$, which are
induced by finitely generated and right-projective idempotent ideals of $A_{n+1}$ and $B_{n+1}$, respectively.
\end{Rem}

{\bf Proof of Theorem \ref{Recollement}.}  We keep all the notation introduced in Lemma \ref{GMCA} and its proof.

$(1)$ By Lemma \ref{LRMOR}, there is a recollement $\big(\D{R_n/I_n}, \D{M_l(R_n,I_n,J_n)}, \D{R_n/J_n}\big)$ induced by a finitely generated, left-projective idempotent ideal of $M_l(R_n,I_n,J_n)$. Thus the recollement restricts to a recollement of bounded-above derived categories. Since $R_n/I_n\simeq A_n\simeq R_n/J_n$ as algebras and since $A_{n+1}$ and $M_l(R_n,I_n,J_n)$ are derived equivalent by Lemma \ref{GMCA}(1), there is a recollement $\big(\Df{A_n}, \Df{A_{n+1}}, \Df{A_n}\big)$.

Similarly, we can apply Lemma \ref{GMCA}(2) and Lemma \ref{LRMOR} to show the existence of the recollement
$\big(\Df{S_n/K_n}, \Df{B_{n+1}}, \Df{S_n/L_n}\big)$. Note that there are isomorphisms of algebras  $S_n/L_n\simeq B_n$ and $$S_n/K_n\simeq (1_{B_n}-f_n)B_n(1_{B_n}-f_n)\simeq (1_{B_{n-1}}-f_{n-1})B_{n-1}(1_{B_{n-1}}-f_{n-1})\simeq\cdots\simeq(1-f_1)B_1(1-f_1)=B_0.$$
This implies the existence of the second recollement in $(1)$.

$(2)$ Note that $R_0$ is symmetric, $A\simeq\End_{R_0}(eA)$ and $D(eA)\simeq Ae$.
Suppose $\dm(A)=\infty$.  By \cite[Lemma 3]{Muller}, $\Ext_{R_0}^i(eA,eA)=0$ for all $i\geq 1$.
It follows from $\Ext_{R_0}^i(eA, eA)\simeq \Ext_{R_0}^i(eA, D(Ae))\simeq D\Tor_i^{R_0}(Ae, eA)$
that $\Tor_i^{R_0}(Ae, eA)=0$ for all $i\geq 1$.
By Proposition \ref{Idem-rec}(3), the recollements in $(3)$ exist for $n=1$.
If $n\geq 1$, then $R_n$ and $S_n$ are symmetric algebras, while $A_n$ and $B_n$ are gendo-symmetric algebras.
Moreover, $\dm(A_n)=\infty=\dm(B_n)$ by $(2)$ and $(1_{B_n}-f_n)B_n(1_{B_n}-f_n)\simeq B_0$ as algebras.
Thus, by induction we can show the existence of recollements for $n\ge 1$. $\square$

\medskip
Theorem \ref{Recollement} can be applied to investigate homological dimensions and higher algebraic $K$-groups.
As usual, for a ring $R$ and $m\in\mathbb{N}$,  we denote by $K_m(R)$ the $m$-th algebraic $K$-group of $R$ in the sense of Quillen, and by $nK_m(R)$ the direct sum of $n$ copies of $K_m(R)$ for $n\ge 0$. If $R$ is an Artin algebra, then $K_0(R)$ is a finitely generated free abelian group of rank $\#(R)$.

\begin{Lem}\label{additivity}
Let $R$ be a ring with $f^2=f\in R$ such that $I:=RfR$ is a strong idempotent ideal of $R$. Suppose that one of the following conditions holds:

$(a)$ Either ${_R}I$ or $I_R$ is finitely generated and projective.

$(b)$ There is a ring homomorphism $\lambda:R/I\to R$ such that the composition of $\lambda$ with the canonical surjection $R\to R/I$ is an isomorphism.

Then $K_n(R)\simeq K_n(fRf)\oplus K_n(R/I)$ for each $n\in\mathbb{N}$.
\end{Lem}

{\it Proof.} When $(a)$ holds, the isomorphisms of algebraic $K$-groups in Lemma \ref{additivity} follow from \cite[Corollary 1.3]{xc4} or \cite[Corollary 1.2]{xc2}.

Let $\pi:R\to R/I$ be the canonical surjection. Clearly, $\pi$ is the universal localization of $R$ at the map $0\to Rf$. Since $I$ is a strong idempotent ideal of $R$, $\pi$ is a homological (also called \emph{stably flat}) ring epimorphism. By \cite[Theorem 0.5]{nr} and \cite[Lemma 2.6]{xc4}, the tensor functors $Rf\otimes_{fRf}-: \prj{(fRf)}\to\prj{R}$ and $(R/I)\otimes_R-:\prj{R}\to\prj{(R/I)}$ induce a long exact sequence of algebraic $K$-groups of rings
$$
\cdots\cdots\to K_{n+1}(R/I)\to K_n(fRf)\to K_n(R)\to K_n(R/I)\to \cdots\to K_0(fRf)\to K_0(R)\to K_0(R/I).
$$
Suppose $(b)$ holds. Then the composition of the functors $R\otimes_{R/I}-: \prj{(R/I)}\to\prj{R}$ with $(R/I)\otimes_R-:\prj{R}\to\prj{(R/I)}$ is an equivalence. This implies that the composition of the maps $K_n(R\otimes_{R/I}-):K_n(R/I)\to K_n(R)$ with $K_n((R/I)\otimes_R-):K_n(R)\to K_n(R/I)$ induced from tensor functors is an isomorphism. Consequently, $0\to K_n(fRf)\to K_n(R)\to K_n(R/I)\to 0$ is split-exact. Thus $K_n(R)\simeq K_n(fRf)\oplus K_n(R/I)$. $\square$

\begin{Koro}\label{FKK}
Let $n$ be a positive integer. Then

$(1)$ $\fd(A_n)\leq \fd(A_{n+1})\leq 2\fd(A_n)+2$ and $\fd(B_0)\leq \fd(B_{n+1})\leq \fd(B_0)+\fd(B_n)+2$.
Thus
$$
\fd(A_{n+1})\leq 2^n\fd(A)+2^{n+1}-2\mbox{ and } \fd(B_{n+1})\leq \fd(A)+n(\fd(B_0)+2).
$$
These inequalities hold true for global dimensions.

$(2)$ $K_*(A_{n+1})\simeq 2^n\,K_*(A)$ and $K_*(B_{n+1})\simeq nK_*(B_0)\oplus K_*(A)$ for $*\in\mathbb{N}$.

$(3)$ If $\dm(A)=\infty$, then
$ K_*(R_n)\simeq K_*(\Lambda)\oplus (2^n-1)K_*(A)$ and $K_*(S_n)\simeq K_*(\Lambda)\oplus nK_*(B_0)$ for any $*\in\mathbb{N}$.
\end{Koro}

{\it Proof.} $(1)$ By Lemma \ref{GMCA}(1), $A_{n+1}$ and $M_l(R_n,I_n,J_n)$ are stably equivalent of Morita type. Since global and finitistic dimensions are invariant under stably equivalences of Morita type, $A_{n+1}$ and $M_l(R_n,I_n,J_n)$ have the same global and finitistic dimensions.
Now, the statements on $A_{n+1}$ in $(1)$ hold by applying \cite[Corollary 3.12 and Theorem 3.17]{xc7} to the recollement $\big(\D{R_n/I_n}, \D{M_l(R_n,I_n,J_n)}, \D{R_n/J_n}\big)$ (see the proof of Theorem \ref{Recollement}(1)). In a similar way, we show the statements on $B_n$ by the recollement $\big(\D{B_0}, \D{B_{n+1}}, \D{B_n}\big)$ in Theorem  \ref{Recollement}(1).

$(2)$ Note that derived equivalent algebras have isomorphic algebraic $K$-groups (see \cite{DS}). By Lemma \ref{additivity} and the proof of Theorem  \ref{Recollement}(1), we have $K_*(A_{n+1})\simeq K_*(M_l(R_n,I_n,J_n))\simeq K_*(R_n/I_n)\oplus K_*(R_n/J_n)\simeq 2K_*(A_n)$ and  $K_*(B_{n+1})\simeq K_*(M_l(S_n,K_n,L_n))\simeq K_*(S_n/K_n)\oplus K_*(S_n/L_n)\simeq K_*(B_0)\oplus K_*(B_n)$. Starting with $A_1=A=B_1$, we can show the isomorphisms in $(2)$ by induction.

$(3)$ By Lemma \ref{additivity} and Theorem \ref{Recollement}(2),   $K_*(R_n)\simeq K_*(R_{n-1})\oplus K_*(A_n)$ and $K_*(S_n)\simeq K_*(S_{n-1})\oplus K_*(B_0)$ for each $n\geq 1$. Together with $(2)$, we can show the isomorphisms in $(3)$ by induction.
$\square$

\begin{Rem}\label{K_0}
Without $\dm(A)=\infty$, the isomorphisms in Corollary \ref{FKK}(3) still hold for $*=0$. This follows from Corollary \ref{FKK}(2) and the fact that if $R$ is a finite-dimensional algebra over a field and $f^2=f\in R$, then $K_0(R)\simeq K_0(fRf)\oplus K_0(R/RfR)$. Thus
$\#(R_n)= \#(\Lambda)+(2^n-1)\,\#(A)$ and $\#(S_n)=\#(\Lambda) + n\; \#(B_0)$.
\end{Rem}

As a consequence of Theorem \ref{Recollement}, we obtain bounds for the stratified dimensions and ratios of  iterated mirror-reflective algebras of gendo-symmetric algebras which are not symmetric.  This provides a new approach to attack the Tachikawa's second conjecture.

\begin{Koro}\label{DMSD}
Let $n$ be a positive integer, and let $(A,e)$ be a gendo-symmetric algebra with $\dm(A)=\infty$. If $A$ is not symmetric, then

$(1)$
$2^n-1\leq \sd(eAe)+(2^n-1)(\sd(A)+1)\leq  \sd(R_n)\leq \#(eAe)+(2^n-1)\,\#(A)-1 \mbox{ and }$
$$n\leq \sd(eAe)+n(\sd(B_0)+1) \leq \sd(S_n)\leq \#(eAe)+ n\,\#(B_0)-1.$$

$(2)$ $\; \frac{\sd(A)+1}{\#(A)}\leq \mathop{\varliminf}\limits_{n\to\infty}
\sr(R_n)\leq 1\; \mbox{ and }\; \frac{\sd(B_0)+1}{\#(B_0)}\leq \mathop{\varliminf}\limits_{n\to\infty}\sr(S_n)\leq 1.$
In particular, if $B_0$ is local, then $\lim\limits_{n\to\infty}\sr(S_n)=1$, where $\varliminf$ means the limit inferior.
\end{Koro}

{\it Proof.} $(1)$ By Theorem \ref{Recollement}(2) and Proposition \ref{upper}(3), $\sd(R_n)\geq \sd(R_{n-1})+\sd(A_n)+1$ and $\sd(S_n)\geq\sd(S_{n-1})+\sd(B_0)+1$.  Similarly, by Remark \ref{Unbounded} and  Proposition \ref{upper}(3), we have $\sd(A_{n+1})\geq 2\, \sd(A_n)+1$, that is, $\sd(A_{n+1})+1\geq 2(\sd(A_n)+1)$.  Moreover, by Proposition \ref{upper}(1), $\sd(R_n)\leq \#(R_n)-1$ and  $\sd(S_n)\leq \#(S_n)-1$. Combining these inequalities with Remark \ref{K_0},  we get $(1)$ by induction.

$(2)$ follows from $(1)$ and Remark \ref{K_0}. $\square$

\medskip
{\bf Proof of Theorem \ref{Tachikawa's second conjecture}.}
$(1)\Rightarrow (2)$  Assume that {\rm (TC2)} holds for all symmetric algebras over $k$. Let  $S$ be an indecomposable symmetric $k$-algebra and $I$ a strong idempotent ideal of $S$.
Then $0=\Ext_{S/I}^i(S/I, S/I)\simeq \Ext_S^i(S/I, S/I)$ for all $i\geq 1$.
This means that ${_S}S/I$ is orthogonal. Then the $S$-module $S/I$ is projective by (1), and therefore ${_S}S\simeq I\oplus S/I$. It follows from $I^2=I$ that $\Hom_S(I,  S/I)=0$.  Since $S$ is symmetric and ${_S}I$ is projective,  $\Hom_S(S/I, I)\simeq D\Hom_S(I,  S/I)=0$. Consequently, $S\simeq \End_S(I)\oplus \End_S(S/I)$ as algebras. Since $S$ is indecomposable, either $\End_S(I)$ or $ \End_S(S/I)$ vanishes. In other words, $I=0$ or $I=S$. This implies that $S$ has no nontrivial strong idempotent ideals. So $(1)$ implies $(2)$.

$(2)\Rightarrow (3)$  An algebra $S$ has no nontrivial strong idempotents if and only if $\sd(S)=0$ if and only if $\sr(S)=0$.  Thus $(3)$ follows.

$(3)\Rightarrow (1)$ Suppose that {\rm (TC2)} does not hold for an indecomposable symmetric algebra $S$ over $k$. Then there exists an indecomposable, non-projective orthogonal $S$-module $M$. Then $A:=\End_S(S\oplus M)$ is a gendo-symmetric, but not a symmetric algebra. Let $S_n$ be the $n$-th reduced mirror symmetric algebra of $A$ for $n\geq 1$. Then $S_n$ is symmetric by Proposition \ref{reduced MS}(1). As $M$ is indecomposable, $\End_S(M)$ is local. Since $M$ contains no nonzero projective direct summands, $S_1$ is indecomposable by Proposition \ref{reduced MS}(3). Further,  by the proof of Theorem \ref{Recollement}(1), $\End_S(M)\simeq (1_{B_n}-f_n)B_n(1_{B_n}-f_n)$ as algebras for any $n\geq 1$. Combining this fact with Proposition \ref{reduced MS}(2), we show that $S_n$ is indecomposable by induction. Since $M$ is orthogonal, we see $\dm(A)=\infty$ by \cite[Lemma 3]{Muller}. It follows from Corollary \ref{DMSD} that  $\lim\limits_{n\to\infty}\sr(S_n)=1$. Thus the supreme in $(3)$ must be $1$, a contradiction to the assumption (3). This shows that $(3)$ implies $(1)$. $\square$

\medskip
{\bf Acknowledgements.} The research was supported partially by the National Natural Science Foundation of China (Grant 12031014, 12122112 and 12171457).

{\footnotesize
\smallskip
Hongxing Chen,

School of Mathematical Sciences  \&  Academy for Multidisciplinary Studies, Capital Normal University, 100048
Beijing, P. R. China;

{\tt Email: chenhx@cnu.edu.cn (H.X. Chen)}

\smallskip
Ming Fang,
Academy of Mathematics and Systems Science, Chinese Academy of Sciences, 100190 \&  School of Mathematical Sciences, University of Chinese Academy of Sciences, 100049
Beijing, P. R. China;

{\tt Email: fming@amss.ac.cn (M. Fang)}

\smallskip
Changchang Xi,

School of Mathematical Sciences, Capital Normal University, 100048
Beijing; \& School of Mathematics and Statistics, Shaanxi Normal University, Xi'an 710062, Shaanxi, P. R. China

{\tt Email: xicc@cnu.edu.cn (C.C. Xi)}

\end{document}